\documentclass[twoside, 12pt, amstex]{amsart}
\makeatletter
\usepackage{graphics}
\usepackage{epsfig}
\usepackage{amsmath}
\usepackage{amssymb}

\begin{document}

\newcommand{\e}{\epsilon}
\newcommand{\ba}{{\mathbf a}}
\newcommand{\bb}{{\mathbf b}}
\newcommand{\bg}{{\mathbf g}}
\newcommand{\bc}{{\mathbf c}}
\newcommand{\w}{{\mathbf w}}
\newcommand{\y}{{\mathbf y}}
\newcommand{\z}{{\mathbf z}}
\newcommand{\x}{{\mathbf x}}
\newcommand{\N}{{\mathbb N}}
\newcommand{\Z}{{\mathbb Z}}
\newcommand{\F}{{\mathbb F}}
\newcommand{\R}{{\mathbb R}}
\newcommand{\Q}{{\mathbb Q}}
\newcommand{\Ql}{{\mathbb Q_{\ell}}}
\newcommand{\C}{{\mathbb C}}
\newcommand{\bL}{{\mathbb L}}
\newcommand{\K}{{\mathbb K}}
\newcommand{\G}{{\mathbf G}}
\newcommand{\BP}{{\mathbb P}}
\newcommand{\cA}{{\mathcal A}}
\newcommand{\cB}{{\mathcal B}}
\newcommand{\cF}{{\mathcal F}}
\newcommand{\cE}{{\mathcal E}}
\newcommand{\cL}{{\mathcal L}}
\newcommand{\cR}{{\mathcal R}}
\newcommand{\cO}{{\mathcal O}}
\newcommand{\cQ}{{\mathcal Q}}
\newcommand{\cX}{{\mathcal X}}
\newcommand{\cH}{{\mathcal H}}
\newcommand{\cM}{{\mathcal M}}
\newcommand{\sB}{{\sf B}}
\newcommand{\cT}{{\mathcal T}}
\newcommand{\cI}{{\mathcal I}}
\newcommand{\cS}{{\mathcal S}}
\newcommand{\sE}{{\sf E}}
\newcommand{\frakp}{{\mathfrak p}}
\newcommand{\frakA}{{\mathfrak A}}
\newcommand{\frakG}{{\mathfrak G}}
\newcommand{\sA}{{\sf A}}
\newcommand{\ga}{{\sf a}}
\newcommand{\es}{{\sf s}}
\newcommand{\m}{{\bold m}}
\newcommand{\bS}{{\bold S}}
\newcommand{\Fp}{\F_p}
\newcommand{\Fq}{\F_q}

\newcommand{\ihra}{\stackrel{i}{\hookrightarrow}}
\newcommand\rank{\mathop{\rm rank}\nolimits}
\newcommand\im{\mathop{\rm Im}\nolimits}
\newcommand\Li{\mathop{\rm Li}\nolimits}
\newcommand\NS{\mathop{\rm NS}\nolimits}
\newcommand\Hom{\mathop{\rm Hom}\nolimits}
\newcommand\Pic{\mathop{\rm Pic}\nolimits}
\newcommand\Spec{\mathop{\rm Spec}\nolimits}
\newcommand\Hilb{\mathop{\rm Hilb}\nolimits}
\newcommand{\length}{\mathop{\rm length}\nolimits}

\newcommand\lra{\longrightarrow}
\newcommand\ra{\rightarrow}
\newcommand\cJ{{\mathcal J}}
\newcommand{\wvskp}{\vspace{1cm}}
\newcommand{\vskp}{\vspace{5mm}}
\newcommand{\nvskp}{\vspace{1mm}}
\newcommand{\nid}{\noindent}
\newcommand{\new}{\nvskp \nid}

\newtheorem{assu}{Assumption}[section]
\newtheorem{obs}{Observation}[section]{
\newtheorem{thm}{Theorem}[section]
\newtheorem{lem}[thm]{Lemma}
\newtheorem{rem}{Remark}[section]
\newtheorem{cor}[thm]{Corollary}
\newtheorem{conj}{Conjecture}[section]
\newtheorem{prop}[thm]{Proposition}
\newtheorem{expl}{Example}[section]
\newtheorem{defn}{Definition}[section]
\newtheorem{quest}{Question}[section]
\newtheorem{prob}{Problem}[section]

\title{Motives and mirror symmetry for Calabi--Yau orbifolds}

\author{Shabnam Kadir}
\address{Institute for Mathematics, University of Hannover, Welfengarten 1,
30167, Hannover, Germany} 
\email{kadir@math.uni-hannover.de}

\author{Noriko Yui}
\thanks{This work is partially supported by N. Yui's
Discovery Grant from the Natural Science Research Council of Canada (NSERC)} 
\address{Department of Mathematics, Queen's University, Kingston
Ontario         Canada K7L 3N6}
\email{ yui@mast.queensu.ca;
yui@fields.utoronto.ca}
\date{June 28, 2006}

\begin{abstract}
We consider certain families of Calabi--Yau orbifolds and their mirror
partners constructed from Fermat hypersurfaces in weighted projective $4$-spaces.
Our focus is the topological mirror symmetry. There are at least three known ingredients 
to describe the topological mirror symmetry, namely,
integral vertices in reflexive polytopes, monomials in graded
polynomial rings (with some group actions), and periods (and Picard--Fuchs
differential equations). In this paper, we will introduce Fermat motives
associated to these Calabi--Yau orbifolds and then use them to
give motivic interpretation of the topological mirror symmetry
phenomenon between mirror pairs of Calabi--Yau orbifolds.
We establish, at the Fermat (the Landau--Ginzburg) point in the moduli space, the
one-to-one correspondence between the monomial classes and Fermat motives. This is done by computing
the number of $\F_q$-rational points on our Calabi--Yau orbifolds
over $\F_q$ in two different ways: Weil's algebraic number theoretic method 
involving Jacobi (Gauss) sums, and Dwork's $p$-adic analytic method involving Dwork characters
and Gauss sums.  We will discuss specific examples in detail.
\end{abstract}

\maketitle
\centerline{Contents}\hfil

\indent\indent 1. Introduction

\indent\indent 2. Fermat motives 

\indent\indent 3. Construction of Calabi--Yau orbifolds 

\indent\indent 4. Construction of mirror Calabi--Yau orbifolds
 
\indent\indent 5. Fermat motives and mirror maps

\indent\indent 6. Batyrev's mirror symmetry 

\indent\indent 7. Monomials and periods

\indent\indent 8. The monomial--motive correspondence : Examples 

\indent\indent 9. Proof of the monomial--motive correspondence

\indent\indent 10. Conclusion and further problems

\indent\indent Bibliography

\vfill
\pagebreak

\section{Introduction}

This is a sequel to the article of Yui [Y05] where Calabi--Yau
orbifolds of Fermat hypersurface threefolds in weighted projective $4$-spaces
were constructed, and their $L$-series (associated to the
$\ell$-adic Galois representations) were determined. It was often the case
that the Galois representations had very high rank, which made it rather
impossible to carry out the required calculations. To remedy this, we introduced
Fermat motives, and then decomposed the Calabi--Yau threefolds into Fermat motives. Via
cohomological realizations of these motives, we were able to 
calculate the motivic $L$-series for each motive. The global $L$-series
was then obtained by gluing the motivic results together.  

These Calabi--Yau orbifolds are all non-rigid ($h^{2,1}>0$), and
their mirror Calabi--Yau threefolds exist satisfying the following
conjecture. 

\begin{conj}[\bf Topological Mirror Symmetry Conjecture]
{\rm Given a family of 
Calabi--Yau threefolds $X$, there is a mirror family of
Calabi--Yau threefolds $\hat X$ such that
$$h^{2,1}(\hat X)=h^{1,1}(X)\quad\mbox{and}\quad h^{1,1}(\hat X)
=h^{2,1}(X)$$
so that the Euler characteristics are subject to the relation
$$\chi(\hat X)=-\chi(X).$$}
\end{conj}

The topological mirror symmetry described above can be reformulated in the
toric geometric setting \`a la Batyrev [Ba94]. 
In Batyrev's theory, mirror symmetry is described in terms of pairs of
reflexive polytopes.  Integral vertices of reflexive polytopes are
the main ingredients in Batyrev's theory, and they correspond to
monomials in graded polynomial rings.  By Aspinwall, Greene and Morrison [AGM93], there is the 
monomial--divisor mirror map (for the corresponding cohomology groups, $H^{1,1}_{toric}$ and 
$H^{2,1}_{poly}$), which yields a one-to-one correspondence between toric divisors of 
a Calabi--Yau family and monomials in the mirror Calabi--Yau family. 

Based on the theory of Dwork, Katz and Griffiths on periods (see Cox and Katz [CK99]),
Candelas et al. [CORV00, CORV03] (resp. Kadir [Ka04, Ka05]) established explicitly a one-to-one 
correspondence between monomials and periods via Picard--Fuchs differential equations,
for the quintic one-parameter (resp. the octic two-parameter) family of Calabi--Yau threefolds.
In their calculations, periods decomposed into the product of subperiods.
This seems to suggest that there should be motivic interpretation for such
factorizations. 

The concept of motives has been emerging in the physics literature, and the purpose
of this paper is to give mathematically rigorous discussion on motives, restricting
our attention to specific examples of Calabi--Yau threefolds. We follow the notion of motives 
due to Grothendieck.  Starting with Fermat hypersurfaces in weighted projective $4$-spaces, 
we define and construct explicitly the so-called Fermat motives from algebraic correspondences,
as described in Shioda [S787]. 
Our goal is to interpret the topological mirror symmetry phenomenon for the mirror pairs of
specific Calabi--Yau orbifolds in terms of Fermat motives and their cohomological
realizations. 
As our main result, we establish a one-to-one correspondence between monomials and
Fermat motives. This correspondence determines Fermat motives which are invariant under the mirror
map.  Since Fermat motives are defined only at the Fermat (the Landau--Ginzburg) point in the 
moduli space, this correspondence is established only at the Fermat point. 

There are more monomials than motives, and we observe that monomials associated to conifold
points seem to be associated to (mixed) Tate motives.  

Incidentally, at the Fermat point, our Calabi--Yau threefolds capture 
the structure of CM type varieties (see Yui [Y05]), and hence our motives are also
CM type motives. 
(G. Moore [Mo98] defined ``attractive'' Calabi--Yau threefolds. Among our
Calabi--Yau orbifolds, there is only one such threefold, namely, 
$m=6, Q=(1,1,1,1,2)$. In fact, this is the only Calabi--Yau orbifold whose weight motive
is rigid, i.e., $h^{2,1}(\cM_Q)=0$ and $B_3(\cM_Q)=2$.)

\smallskip

Now we will describe the contents of this paper.

Section \ref{Fermatmotives} is devoted to the definition of Fermat motives and their 
cohomological realizations. We use the definition of motives due to Grothendieck, 
and Manin [Ma70], which is based on algebraic correspondences and projectors. We follow the 
exposition of Shioda [Sh87] and Gouv\^ea and Yui [GY95]. 

In Section \ref{constructionorbifolds}, we construct Calabi--Yau orbifolds in weighted 
projective $4$-spaces. The starting point is the Fermat hypersurface $V$ of
degree $m\geq 5$ and dimension $3$, and a finite abelian groups (which is  
a subgroup of the automorphism group of $V$). This group will determine a weight.
We take the quotient of $V$ by such a group. This gives rise to a quotient threefold with 
singularities.  We then resolve singularities by taking the crepant resolution (which is 
guaranteed to exist for dimension $\leq 3$). The smooth threefold thus obtained is our   
Calabi--Yau threefold. There are altogether $147$ such Calabi--Yau orbifolds.

In Section \ref{constructionmirrororbifolds}, we describe the Greene--Plesser orbifolding 
construction of mirror partners of the Calabi--Yau orbifolds constructed in the previous section. 
We review the mirror construction from the paper of Greene and Plesser [GP90]. 
The mirror symmetry is interpreted as the duality between the two finite abelian groups
associated to the mirror pair of Calabi--Yau threefolds.

In Section \ref{Fermatmotivesandmirrormaps}, we will construct Fermat motives for the 
Calabi--Yau orbifolds in Section \ref{constructionorbifolds}, and compute their invariants 
(e.g., motivic Hodge numbers, motivic Betti numbers) via their cohomological realizations.  For each 
mirror pairs of Calabi--Yau threefolds, we also determine Fermat motives which are 
invariant under the  mirror maps.  In particular, we observe that for each Calabi--Yau orbifold, 
the motive associated to the weight is always invariant under the mirror map.
If $h^{1,1}=1$, the weight motive is the only motive invariant under the mirror map.
However, when $h^{1,1}>1$, there are other motives apart from the weight motive that
remain invariant under the mirror map.  

In Section \ref{Batyrevmirrorsymmetry}, we review the construction of mirror pairs 
of Calabi--Yau hypersurfaces in toric geometry due to Batyrev [Ba94]. We will confine ourselves 
to Calabi--Yau threefolds. 
Reflexive polytopes and their dual polytopes are the main players in Batyrev's
toric mirror symmetry. A pair of reflexive polytopes $(\Delta,\Delta^*)$ 
gives rise to a mirror pairs of Calabi--Yau hypersurfaces. 
It is noted that the origin is the only integral point contained both in the reflexive
polytope and its dual polytope. (This fact plays a pivotal role in proving our
main result.)  Integral points correspond to monomials
in graded polynomial rings. We will discuss, in particular, the monomial--divisor 
mirror map of Aspinwall, Greene and Morrison [AGM93], which gives the isomorphism between 
the two spaces $H_{toric}^{1,1}(X)$ and $H^{2,1}_{poly}(\hat X)$
for a mirror pair $(X, \hat X)$. This establishes a one-to-one correspondence between
integral points in the reflexive polytope of $X$ and monomials in the polynomial ring of $\hat X$.  

In Section \ref{monomialsandperiods}, we will discuss a one-to-one correspondence between monomials
and periods via Picard--Fuchs differential equations.  The method
of Dwork--Katz--Griffiths determines the Picard--Fuchs differential equations
for Calabi--Yau hypersurfaces in weighted projective spaces (of any dimension) (see
Cox and Katz [CK99]). In this section, we will illustrate this
correspondence focusing on the concrete calculations of periods of Candelas et al. 
[CORV00, CORV03] for the one-parameter deformation of the quintic Calabi--Yau threefold 
in the ordinary projective $4$-space $\BP^4$, and of Kadia [Ka05] for the two-parameter 
deformation of the octic Calabi--Yau threefold in the weighted projective space $\BP^4(Q)$ 
with $Q=(1,1,2,2,2)$.  Here we observe that the Picard--Fuchs differential quation decomposes
into the product of smaller order Picard--Fuchs differential equations. It suggests that such a 
decomposition ought to have origin in the motivic decomposition of the manifold.

The Section \ref{monomialFermatcorrespondence} contains our main examples and the main result
on the monomial--motive correspondence (Theorem 8.1). 
We establish a one-to-one correspondence between the class of monomials and Fermat motives at 
the Fermat point for the Calabi--Yau threefolds of Section \ref{constructionorbifolds}.
We prove that the motives which are invariant under the mirror maps correspond to
the class of the constant monomial (and hence to the origin in the polytopes). 
We illustrate the monomial--motive correspondence for the quintic and the octic
Calabi--Yau orbifolds.

The Section \ref{proofcorrespondence} contains our proof for the monomial--motive 
correspondence.
We compute the number of $\F_q$--rational points (and hence congruence zeta-functions) 
for our Calabi--Yau threefolds over finite
fields $\F_q$ in two different ways. On one hand, we compute them with
Weil's method using Jacobi (Gauss) sums. On the other hand, we compute them with
Dwork's $p$-adic method using Dwork's characters and Gauss sums. We then
show that the two approaches reconcile at the Fermat point.  

Away from the Fermat point, Calabi--Yau orbifolds with deformation parameters
yield more monomials than motives. A conifold point on these Calabi--Yau orbifolds with deformation 
parameters is a singularity locally isomorphic
to the projective quadric surface $X^2+Y^2+Z^2+T^2=0$, and we may associate a motive 
(e.g., a mixed Tate motive) employing the same line of arguments as in Bloch--Esnault--Kreimer 
[BEK05].

Finally, Section \ref{conclusions} presents conclusions of this work, and further
problems and future projects.  The main conclusion is to bring in motives to
the realm of the topological mirror symmetry for our Calabi--Yau orbifolds. 
\medskip

\hfil {\bf Acknowledgments}\hfil
\smallskip

Shabnam Kadir held a Fields Postdoctoral Fellowship in 2004--2005 supported by the Fields
Institute Thematic Program ``Geometry of String Theory'', and N. Yui's
NSERC Discovery Grant.  

During the course of this work, Noriko Yui held visiting positions at 
a number of institutions.  This include University of  
Leiden (Winter 2004), Max-Planck-Institut-f\"ur-Mathematik 
Bonn (Spring 2004, Spring 2006), University of Hannover (Spring 2004, 2006), 
University of Mainz (Spring 2004),
Tsuda College (Spring 2005) and Kavli Institute for Theoretical
Physics (KITP) (December 2005) at University of California Santa Barbara. The research at KITP 
was supported in part by the National Science Foundation under the Grant No. PHY99-07949.
She is grateful for the hospitality of these institutions. 
\medskip

\section{Fermat motives}\label{Fermatmotives}

We will employ the definition of motives due to Grothendieck and Manin [Ma70]. 
The most recommended references for the generality of motives might be Deligne--Milne [DM82],
and Soul\'e [So84].  In this paper we will be confining ourselves to the so-called
``Fermat'' motives arising from Fermat hypersurfaces. We will recall the construction of 
Fermat motives from Shioda [Sh87] and Gouv\^ea and Yui [GY95].  The construction works for
any dimension, not only for dimension $3$. 

We start with the Fermat hypersurface of degree $m$ and of dimension $n$ in the
projective space $\BP^{n+1}$:
$$V\,:\, Z_1^m+Z_2^m+\cdots+Z_{n+2}^m=0\in\BP^{n+1}$$
Let $\mu_m$ be the group of $m$-th roots of unity and let
$$\mathfrak{G}:=(\mu_m)^{n+2}/(\mbox{diagonal})=
\{\bg:=(g_1,g_2,\cdots, g_{n+2})\,|\, g_i\in \mu_m\quad\forall i\}/\{(g,g,\cdots,g)\}$$
be a subgroup of the automorphism group $\mbox{Aut}(V)$ of $V$. The group $\mathfrak G$ is
a finite group of order $m^{n+1}$, and it acts on $V$ by component-wise multiplication of $g_i$ for each $i$.   
The character group of $\mathfrak{G}$ is identified with the set  
$$\hat{\mathfrak{G}}:=
\{\ba=(a_1,a_2,\cdots, a_{n+2})\,|\, a_i\in (\Z/m\Z)\,\,\forall\, i,\, \sum_{i=1}^{n+2} a_i \equiv 0 \pmod m\,\}$$
and the duality between $\mathfrak{G}$ and $\hat{\mathfrak{G}}$ is
given by:
$$\hat{\mathfrak{G}}\times \mathfrak{G}\mapsto \bL=:\Q(\zeta_m)\,\quad\,
(\ba,\, \bg)\mapsto \ba(\bg):=\prod_{i=1}^{n+2} g_i^{a_i}$$
where $\zeta_m=e^{2\pi i/m}$ is a primitive $m$-th root of unity.
For each $\ba\in\hat{\mathfrak{G}}$, let $A=[\ba]$ be the $(\Z/m\Z)^{\times}$-orbit of $\ba$.
Put $\bL_{\ba}=\bL_A=:\Q(\zeta_m^{d})$ where $d=:\mbox{gcd}(\ba, m)$. 

\begin{defn}\label{projector}
{\rm For each character $\ba\in\hat{\mathfrak{G}}$, define
$$p_{\ba}:=\frac{1}{\#{\hat{\mathfrak{G}}}}\sum_{\bg\in\mathfrak{G}}\ba(\bg)^{-1}\bg
=\frac{1}{m^{n+1}}\sum_{\ba\in\mathfrak{G}}\ba(\bg)^{-1}\bg,$$
and for each $(\Z/m\Z)^{\times}$-orbit $A=[\ba]$ of $\ba$, define
$$p_A:=\sum _{\ba\in A} p_{\ba}
=\frac{1}{m^{n+1}}\sum_{\bg\in\mathfrak{G}}(\mbox{trace}_{\bL_{\ba}/\bL}(\ba(\bg)^{-1}))\,\bg.$$
Then it is easy to see that $p_{\ba}$ and $p_A$ are elements of the group ring $\bL[\mathfrak{G}]$
and $\Z[\frac{1}{m}][\mathfrak{G}]$, respectively.}
\end{defn}

\begin{prop}
$p_{\ba}$ and $p_A$ are projectors (idempotents), that is,
$$p_{\ba}\cdot p_{\bb}=\begin{cases} p_{\ba}&\quad\mbox{if $\ba=\bb$} \\
                                        0   &\quad\mbox{if $\ba\neq \bb$}\end{cases}
$$
and
$$p_A\cdot p_B=\begin{cases} p_A &\quad\mbox{if $A=B$}\\
                              0 & \quad\mbox{if $A\neq B$}\end{cases}$$
Furthermore, we have the decomposition
$$\sum_{\ba\in{\hat{\mathfrak{G}}}} p_{\ba}=1,\quad\mbox{and}\quad
\sum_{A\in O(\hat{\mathfrak{G}})} p_A=1$$
where $O(\hat{\mathfrak{G}})$ denotes the set of $(\Z/m\Z)^{\times}$-orbits in $\hat{\mathfrak{G}}$.
\end{prop}
 
\begin{defn} {\rm  Identifying each automorphism $\bg\in\mbox{Aut}(V)$ with its graph, 
the projector 
$p_A$ may be regarded as an algebraic $n$-cycle on $V\times V$ with coefficients in
$\Z[\frac{1}{m}]$.  The pair $$(V,\, p_A)=:\cM_A$$ may be called the {\it Fermat motive}
corresponding to the $(\Z/m\Z)^{\times}$-orbit $A$ of $\hat{\mathfrak{G}}$. 
Furthermore, there is the motivic decomposition of the Fermat hypersurface:
$$(V, {\bf 1})=\oplus_{A\in O(\hat{\mathfrak{G}})} (V, p_A)=\oplus_{A\in O(\hat{\mathfrak{G}})}\cM_A
$$
corresponding to $\sum_{A\in O(\hat{\mathfrak{G}})} p_A=1$. (Here ${\bf 1}:=\Delta_V$.)

Note that Fermat motives are well-defined over $\Q$ or finite fields $\F_p$ for prime $p$
such that $(m,p)=1$.}
\end{defn} 

Now we will discuss the cohomological realizations of Fermat motives. 

\begin{defn} {\rm Let $H^{\bullet}(V,\star)$ be any Weil cohomology group. For a Fermat motive
$\cM_A=(V, p_A)$, define
$$H^{\bullet}(\cM_A,\star):=H^{\bullet}(V, \star)^{p_A}$$
as the image of $p_A$ (or equivalently, the kernel of $p_A-1$) acting on $H^{\bullet}(V,\star)$.}
\end{defn}

The motivic decomposition yields the decomposition of cohomology groups
$H^{\bullet}(V,\star)$ with various coefficients $\star$, e.g., the $\ell$-adic \'etale cohomology,
the de Rham cohomology, the crystalline cohomology, etc.  
Here we will discuss de Rham and Betti realizations.

\begin{prop}
Let
$$\mathfrak{A}:=\{\ba=(a_1,a_2,\cdots,a_{n+2})\in\hat{\mathfrak{G}}\,|\, a_i\neq 0\in (\Z/m\Z) \,\forall\,i\,\}$$
be a subset of $\hat{\mathfrak{G}}$.
For $\ba=(a_1,a_2,\cdots, a_{n+2})\in\mathfrak{A}$, let
$$\Vert \ba \Vert:=\left(\frac{1}{m}\sum_{\ba\in A} a_i\right)-1.$$
 
(a) The $(i,n-i)$-th {\it motivic Hodge number} is given by 
$$h^{i,n-i}(\cM_A)=\begin{cases}
\#\{\ba\in A\,|\, \Vert\ba\Vert=i& \quad \mbox{for $i,\,0\leq i\leq n$ and $A\subset\mathfrak{A}$}\\
                                1& \quad \mbox{for $i=\frac{n}{2}$ and $A=[(0,0,\cdots,0)]$}\\
                                 0&\quad \mbox{otherwise}\end{cases}
$$
In particular, the {\it geometric genus} of $\cM_A$ is
$p_g(\cM_A):=h^{0,n}(\cM_A)$.
We have the Hodge decomposition
$$h^{i,n-1}(V)=\sum_A h^{i,n-i}(\cM_A)$$
and in particular,
$$p_g(V)=\sum_{A} p_g(\cM_A).$$

(b) The $i$-th {\it motivic Betti number} is given by
$$B_i(\cM_A)=\begin{cases} \#A & \quad\mbox{for $i=n$ and $A\subset \mathfrak{A}$}\\
                           1 & \quad\mbox{for $i$ even and $A=[(0,0,\dots,0)]$} \\
                           0 &\quad\mbox{otherwise}\end{cases}
$$
For $i=n$, we have 
$$B_n(\cM_A)\leq \varphi(m)/\mbox{gcd}(\ba, m)$$
where $\varphi(m)$ denotes the Euler's phi-function. (We may call $B_n(\cM_A)$
the {\it dimension} of the motive $\cM_A$.)
 
The global $i$-th Betti number $B_i(V)$ is given by   
$$B_i(V)=\sum_A B_i(\cM_A).$$

For $h^{i,n-i}(V)$ and $B_i(V)$, the sum runs over $A$ in the set $O(\mathfrak{A})$ of $(\Z/m\Z)^{\times}$-orbits in
$\mathfrak{A}$ if $n$ is odd, and the set $O(\mathfrak{A}\cup\{0\})$ if $n$ is even.
\end{prop}

\section{Construction of Calabi--Yau orbifolds}\label{constructionorbifolds}

From now on, we will confine ourselves to threefolds ($n=3$).
We will construct mirror pairs of Calabi--Yau threefolds by orbifolding method due to
Greene and Plesser [GP90].

Fix a Fermat hypersurface threefold $V=V(m,3,{\bf 1})$ of degree $m$ and coefficients $\bf 1$ defined over $\Q$:
$$V: Z_1^m+Z_2^m+Z_3^m+Z_4^m+Z_5^m=0\subset\BP^4_{\Q}.$$
Let $Q=(q_1,q_2,q_3,q_4,q_5)\in{\N}^5$ be a vector consisting 
of $5$-tuples of positive rational integers.  We call $Q$ a {\it weight}
if it satisfies the following two conditions: 
\smallskip

(1) ${\mbox{gcd}}(q_1,\cdots, \widehat{q_j},\cdots, q_5)=1$  
for every $j,\, 1\leq j\leq 5$, and 

(2) $q_j\,|\,m=\rm{deg}(V)$ for every $j,\,1\leq j\leq 5$. 

(Here $\hat{\quad}$ means the exclusion of that component.)
\smallskip

For each $j,\,1\leq j\leq 5$, let $\mu_{q_j}$ denote the finite
group scheme of $q_j$-th roots of unity over $\Q$, namely,
$$\mu_{q_j}:={\mbox{Spec}}\,(\Q[T]/(T^{q_j}-1)).$$ 
Taking the direct product of $\mu_{q_j}$'s and define
$$ G=G_Q:=\mu_{q_1}\times \cdots\times \mu_{q_5}.$$
Then $G=G_Q$ is a finite automorphism group scheme of $V_{\Q}$, whose
action over the closure $\bar{\Q}$ is described as follows:
$$
(Z_1,\cdots, Z_5)\mapsto (\zeta_{q_1}^{e_1}Z_1,\cdots, 
\zeta_{q_5}^{e_5}Z_5), \quad \text{$e_j\in\Z/q_j\Z$ for each $j$}$$ 
where $\zeta_{q_j}$ denotes a primitive $q_j$--th root of unity in 
$\bar{\Q}$ for every $j,\,1\leq j\leq 5$. 

Take the quotient, $Y:=V/{G_Q}$.  Then $Y$ is, in general, a singular variety 
defined over $\Q$ in the weighed projective $4$-space $\BP^4(Q)$. Let $\Sigma(Y)$ denote
the singular loci of $Y$. Then $\Sigma(Y)$ is finite, and all singularities are at most abelian (indeed, cyclic) 
quotient singularities. 

Next we wish to construct a smooth Calabi--Yau threefold from $Y$. 
For this we first look for the sufficient condition for which this quotient variety $Y$ is a singular 
Calabi--Yau threefold, that is, the canonical sheaf $\omega_Y$ to be trivial.
The so-called Calabi--Yau condition is imposed on the weight $Q=(q_1,q_2,\cdots, q_5)$ by 
\smallskip

(3) $q_1+q_2+q_3+q_4+q_5=m\quad\mbox{or equivalently}\quad \sum_{i=1}^5\frac{q_i}{m}=1.$ 
\smallskip

Indeed, in the weighted projective $4$-space $\BP^4(Q)$, the canonical sheaf is given by (see
Dolgachev [Dol82])
$$\omega_Y\simeq\mathcal{O}(m-q_1-q_2-q_3-q_4-q_5).$$
We require the triviality of $\omega_Y$, which gives rise to the sufficient condition (3).
With the condition (3), $Y$ has trivial canonical bundle away from its singularities. 

\begin{defn} {\rm
For a fixed $m$, we call $Q=(q_1,q_2,q_3,q_4,q_5)$ an {\it admissible} weight if it satisfies
the conditions (1), (2) and (3), and we will call $<m, Q>$ an {\it admissible} pair.}
\end{defn}

Now we look for a smooth resolution of singularities without disturbing the triviality of the canonical bundle.
Such a resolution is called a {\it crepant} resolution. The existence of the crepant resolution for our
singular Calabi--Yau threefold $Y$ is proved by Greene, Roan and Yau [GRY91]. 
We summarize their results reformulated by Yui [Y05] suitable for the arithmetic discussion
and applications.

\begin{thm}\label{thm3.1} {\sl Let $V=V(m,3,{\bf 1})$ be the Fermat hypersurface
threefold of degree $m\geq 5$ with coefficients $\bf 1$ defined over $\Q$.  
Let $Q=(q_1,q_2,q_3,q_4,q_5)$ be a weight.

(a) There are $147$ admissible pairs $<m, Q>$ giving rise to singular
Calabi--Yau threefolds $Y$. The smallest is $<5, (1,1,1,1,1)>$ and the largest
is $<1807, (1, 42, 258, 602, 903)>$.

In the weighted projective $4$-space $\BP_{\Q}^4(Q)$, $Y$ is
defined by the equation
$$Y_1^{m/q_1}+Y_2^{m/q_2}+Y_3^{m/q_3}+Y_4^{m/q_4}+Y_5^{m/q_5}=0
\subset\BP_{\Q}^4(Q)$$
with ${\rm{deg}}(Y_j)=q_j$ for each $j,\,1\leq j\leq 5$.

(b) $Y$ has at most abelian (cyclic) quotient singularities. Let $\Sigma(Y)$
be the singular loci of $Y$. 
Let $\ell$ be an integer ranging over the set $\{0,1,\cdots, m-1\}$.
For each $\ell$, we define the two sets of indices: 
$$S_{\ell}:=\{\,j\,| \,1 \leq j\leq 5, \ell/m_j\in\Z\,\}\quad
\rm{and}\quad I_{\ell}=\{1,2,3,4,5\}\setminus S_{\ell}.$$
Further put $c_{I_{\ell}}:={\mbox{gcd}}\{q_j\,|\, j\not\in I_{\ell}\}
={\mbox{gcd}}\{q_j\,|\,j\in S_{\ell}\}$.  
Let
$$Y_{I_{\ell}}:=Y\cap \{Y_j=0,\, j\in I_{\ell}\}.$$
Then $Y_{I_{\ell}}\in\Sigma(Y)$, and it is defined over $\Q$.
\smallskip

More concretely, the following assertions hold: 
\smallskip

(b1) $\# S_{\ell}\neq 4$ for any $\ell$, and $\#S_{\ell}=5$ for
$\ell=0$.

(b2) If $\# S_{\ell}=3$ for some $\ell$, then there is a set $I_{\ell}$ with
$\#I_{\ell}=2$, $c_{I_{\ell}}\geq 2$. In this case $\Sigma(Y)$ contains 
a dimension $1$ singularity. 

Suppose that $C$ is an irreducible component of 
$Y_{I_{\ell}} =Y\cap\{Y_j=0\,|\,j\in I_{\ell}\},$ then $C$ is a smooth weighted 
diagonal curve defined over $\Q$ of degree 
$m^{\prime}:=m/c_{I_{\ell}}$ and reduced weight 
$Q^{\prime}:=(q^{\prime}_i,q^{\prime}_j,q^{\prime}_k)$ where 
$i,j,k\not\in I_{\ell}$ and $q^{\prime}_{\bullet}=q_{\bullet}/c_{I_{\ell}}$.  
The multiplicity of $C$ in $Y_{I_{\ell}}$ is $m_C:=c_{I_{\ell}}-1$.
The genus $g(C)$ of $C$ is given by the coefficient of 
$t^{m^{\prime}-(q^{\prime}_i+q^{\prime}_j+q^{\prime}_k)}$ in the formal 
power series 
$$\frac{1-t^{m^{\prime}}}{(1-t^{q^{\prime}_i})(1-t^{q^{\prime}_j})
(1-t^{q^{\prime}_k})}.$$  

\smallskip

(b3) If $\# S_{\ell}=2$ for some $\ell$, then there are dimension $0$
singularities.  Let $P\in\Sigma(Y)$ be a singular point. Then
$P$ is a cyclic quotient singularity, and $\pi^{-1}(P)$ is
a smooth rational surface birationally equivalent to $\BP^2$. 
It is defined over a finite Galois extension of $\Q$; its Galois
orbit is defined over $\Q$.
\smallskip

(b4) If $\# S_{\ell}=1$, $Y$ has no singularity. 
\medskip

(c) There exists the crepant resolution $\pi\,:\, X\to Y$, and $X$ is a smooth
Calabi--Yau threefold defined over $\Q$.
The third Betti number of $X$ can be computed by
$$B_3(X)=B_3(Y)+\sum_C m_C\, B_1(C)$$
where $C$ runs over all the smooth irreducible curves in $\Sigma(Y)$ with 
multiplicity $m_C=c_{I_{\ell}}-1$.
($\pi^{-1}(C)$ is a smooth ruled surface birationally equivalent to $C\times\BP^1$
and it is defined over $\Q$.)
\smallskip

Furthermore, $B_3(Y)$ can be computed by Vafa's formula (see Roan [Ro90]) as follows.
For each integer $\ell\in \{0,1,2,\cdots, m-1\}$, 
let $S_{\ell}=\{\,j\,\vert 1\leq j\leq 5, \ell\frac{q_i}{m}\in\Z\}$ be the
set of indices defined in Proposition 3.1, 
and define the integer $\beta_{\ell}$ by
$$\beta_{\ell}:=\frac{1}{m}\sum_{r=0}^{m-1} \prod_{\ell\frac{q_i}{m}\in\Z, 
r\frac{q_i}{m}\in\Z} \left(1-\frac{m}{q_i}\right).$$
Here we employ the convention that the product 
$\prod_{\ell\frac{q_i}{m}, r\frac{q_i}{m}} (1-\frac{m}{q_i})=1$ if there is no  
$q_i$ with $\ell\frac{q_i}{m}\in\Z, r\frac{q_i}{m}\in\Z$. 
Then
$$B_3(Y)=-\beta_0,\quad{\mbox{and}}\quad 
B_3(X)=B_3(Y)-\sum_{\# s_{\ell}=3} \beta_{\ell},$$
where the sum runs over the set $s_{\ell}$ consisting of three
elements.  The Euler characteristic of $X$ is given by 
$$\chi(X)=\sum_{\ell=0}^{m-1} \beta_{\ell}.$$
Consequently,
$$B_2(X)=h^{1,1}(X)=
-1+\frac{1}{2}\sum_{\ell=0}^{m-1}\sum_{\# s_{\ell}<3} \beta_{\ell}.$$
In particular, all cycles generating $H^{1,1}(X)$ arise from singular
points on $Y$.}
\end{thm}

\begin{rem} {\rm Here we used  Vafa's method developed for computing 
orbifold Euler characteristics. (At glance, it is not transparent why $\beta_{\ell}$ is
an integer. Hendrik Lenstra gave an elementary proof to establish the
integrality of $\beta_{\ell}$'s.)  It should be noted that Vafa's method is
a special case of the calculation of {\it stringy} Euler characteristics by toric
method developed by Batyrev [Ba94]. 
In Section 5 we will give another approach for the computation of Betti
numbers and Euler characteristics of the 147 Calabi--Yau orbifolds using Fermat motives.}
\end{rem}

\section{Construction of mirror Calabi--Yau orbifolds}\label{constructionmirrororbifolds}

The Calabi--Yau orbifolds constructed in Theorem 3.1 are not rigid,
i.e., $h^{1,2}\neq 0$. Therefore, their mirror Calabi--Yau threefolds exist in the sense of Topological
Mirror Symmetry Conjecture 1.1.
In this section, we will construct mirror partners of these Calabi--Yau
orbifolds by applying again orbifolding construction on perturbed Calabi--Yau hypersurfaces by certain 
(finite abelian) groups of discrete symmetries.  This construction is due
to Greene and Plesser [GP90] (see also Greene--Roan and Yau [GRY91], Roan [Ro90], [Ro91], [Ro94]). 
Now we will describe the mirror construction.

\begin{thm}\label{thm4.1} {\sl
Let $X$ be a Calabi--Yau orbifold corresponding to $<m,Q>$ with $Q=(q_1,q_2,q_3,q_4,q_5)$ in Theorem 3.1 
obtained from the weighted Fermat hypersurface:
$$Y_1^{m_1}+Y_2^{m_2}+Y_3^{m_3}+Y_4^{m_4}+Y_5^{m_5}=0 \subset{\BP}^4(Q)$$
where we put $m_j=m/q_j$ for $1\leq j\leq 5$.  Let ${\mathcal{W}}_Q$ be the generic weighted
Fermat hypersurface with deformation parameters, and let $\{{\mathcal{W}}_Q=0\}$ be its zero locus in ${\BP}^4(Q)$.
Then $\{{\mathcal{W}}_Q=0\}$ is a Calabi--Yau threefold.
Let $\Omega$ denote the unique (up to a scalar multiplication) holomorphic 
$3$-form on ${\mathcal{W}}_Q$.  Then
$$\Omega=\mbox{Res}_{{\mathcal{W}}_Q=0} d\frac{d\mu}{{\mathcal{W}}_Q}$$
where 
$$d\mu=\sum_{j=1}^5 (-1)^j q_jY_jdY_1\cdots\wedge\widehat{dY_j}\wedge\cdots 
\wedge dY_5.$$

For each $j,\, 1\leq j\leq 5$, let $g_{m_j}$ be a primitive $m_j$-th 
root of unity and let $\mu_{m_j}$ be the cyclic group generated by 
$g_{j}=e^{2 \pi i/m_j}$.   
Define the group $\hat G$ by
$$\hat G=\{\,\bg=(g_1,g_2,g_3,g_4,g_5)\,|\, g_i^{m_i}=1,\, \prod_{i=1}^5 g_i=1\,\}/
\{g,g,\cdots,g)\}.$$
Then $\hat G$ is a finite abelian group (scheme) of order 
$\prod_{j=1}^5 m_j/m^2=m^3/\prod_{i=1}^5 q_i$, and it leaves $\Omega$ invariant. 

This gives rise to a $\hat G$-invariant hypersurface
$${\mathcal{M}}_Q:=\sum_{I=(i_1,\cdots, i_5)} c_I Y_1^{i_1}\cdots Y_5^{i_5}$$ 
in ${\BP}^4(Q)$ quotiented out by $\hat G$, 
where the sum runs over $$I=(i_1,\cdots, i_5)\in(\Z/m_1\Z)\times
\cdots \times (\Z/m_5\Z)$$ subject to the relation $q_1i_1+\cdots+q_5i_5=m$.
The zero locus $\{{\mathcal{M}}_Q=0\}$ defines a Calabi--Yau threefold. 

Then there exists a resolution  
$\hat X:=\{{\mathcal{M}}_Q=0\}\subset{\BP}^4(Q)/\hat G$ which is a mirror partner
of the family of Calabi--Yau orbifold $X:=\{{\mathcal{W}}_Q=0\}$ satisfying
Topological Mirror Symmetry Conjecture 1.1: 
$$h^{2,1}(\hat X)=h^{1,1}(X),\, h^{1,1}(\hat X)=h^{2,1}(X),\, 
\chi(\hat X)=-\chi(X).$$} 
\end{thm}

The orbifolding construction was described by Greene and
Plesser [GP90]. The calculations of Hodge numbers and Euler characteristics
were carried out in Candelas-Lynker-Schimmrigk [CLS90]. 
Roan [R091] gave a mathematical proof of these results.

\begin{rem}
{\rm In Theorem 4.1, defining equations for mirror Calabi--Yau
threefolds are given by deformations of weighted Fermat
hypersurfaces threefolds. 
Generically, these defining equations
may have as many independent deformation parameters as 
$h^{2,1}_{poly}< h^{2,1}$. (The Hodge number $h^{2,1}_{poly}$ is defined
in terms of toric geometry. For the definition, see Sections 6 and 7 below.)
The Hodge numbers are invariant under K\"ahler deformations, i.e, when deformations
preserve K\"ahler structure of the manifold in question (see
Nakamura [Na75]), so $h^{1,1}$ remains unchanged under deformations. 
In particular, by taking a special point of the generic deformation,
we can have a defining equation for our Calabi--Yau family with a small number of
deformation parameters.}
\end{rem}  

\begin{expl}{\rm
We consider one-parameter deformations of 
Calabi--Yau threefolds obtained from weighted Fermat hypersurfaces 
corresponding to a pair $<m,Q>$ in Theorem 4.1:
$$Y_1^{m/q_1}+Y_2^{m/q_2}+Y_3^{m/q_3}+Y_4^{m/q_4}
+Y_5^{m/q_5}-m\,\psi Y_1Y_2Y_3Y_4Y_5=0\subset{\BP}^4(Q)$$
where $\lambda$ is a parameter (subject to the relation $\psi^m=1$).  
Then its zero locus ${\mathcal{W}}_Q=0$ gives rise to a
family of Calabi--Yau threefolds. 

TABLE 1 lists Calabi--Yau threefolds $\{{\mathcal{W}}_Q=0\}$ 
corresponding to the pairs $<m,Q>$ with $h^{1,1}(X)=1$,
and their mirror families of Calabi--Yau threefolds
${\mathcal{M}}_Q$. Here a generator $\bg=(g_1,\cdots, g_5)$ for
$\hat G$ stands for $(e^{2 \pi i g_1/m},\cdots, e^{2 \pi ig_5/m})$. 
The action of $\hat G$ on $X$ is
$$(Y_1,Y_2,Y_3,Y_4,Y_5)\mapsto (g_1Y_1,g_2Y_2,g_3Y_3,g_4Y_4,g_5Y_5).$$
For instance, for $m=5$ and the generator $(1,0,0,4,0)$, the
action is read as $$(Y_1,Y_2,Y_3,Y_4,Y_5)\mapsto
(e^{2\pi i/5}Y_1, Y_2,Y_3,e^{2\pi i 4/5}Y_4,Y_5).$$} 
\end{expl}

\begin{table}
\begin{center}
\begin{tabular}{|c|c||c|c|c|} \hline
$m$ & $Q$ & $G_Q$ & $\hat G$ &  $\mbox{generators}$ \\ \hline 
$5$ & $(1,1,1,1,1)$ & $\{1\}$  & $(\Z/5\Z)^3$ & $(1,0,0,4,0),(1,0,4,0,0),(1,4,0,0,0)$
\\ \hline
$6$ & $(1,1,1,1,2)$ & $(\Z/2\Z)$ & $(\Z/3\Z)\times(\Z/6\Z)^2$ & $(0,2,2,2,0),(5,0,0,1,0),
(0,5,0,1,0)$ \\ \hline 
$8$ & $(1,1,1,1,4)$ & $(\Z/4\Z)$ & $(\Z/2\Z)\times(\Z/8\Z)^2$ & $(7,0,0,1,0),
(0,7,0,1,0),(0,0,7,1,0)$ \\ \hline
$10$ & $(1,1,1,2,5)$ & $(\Z/2\Z)\times(\Z/5\Z)$ & $(\Z/10\Z)^2$ &$(9,0,1,0,0,),(0,9,1,0,0)$ \\ \hline
\end{tabular}
\vskip 0.2cm
\caption{Mirror Calabi--Yau orbifolds I}
\end{center}
\end{table}

The above construction of mirror pairs of Calabi--Yau orbifolds may be characterized in terms of 
the ``duality'' between the two finite abelian groups. See Greene--Plesser [GP90],
and for a nice exposition, Morrison [Mor97].

\begin{thm}
{\sl For the above mirror pairs of Calabi--Yau orbifolds $(X, \hat X)$, 
the mirror symmetry can be described 
as the ``duality'' between the two finite abelian groups, that is,
there is the group $G=G_Q$ associated to the original Calabi--Yau orbifold $X$:
$$G=G_Q=\prod_{j=1}^5 \mu_{q_j}\quad\mbox{of order\quad $\prod_{j=1}^5 q_j$}.$$
On the mirror side, there is the group $\hat G$ associated to the mirror
Calabi--Yau orbifold $\hat X$: 
$$\hat G=\{\bg=(g_1,\cdots, g_5)\,|\, g_j^{m_j}=1,\, 
\prod_{j=1}^5 g_j=1\}/\mu_m\quad\mbox{of order\quad $m^3/\prod_{j=1}^5 q_j$}.$$
These two abelian groups are subject to the ``duality'' relation:
$$\#G_Q\times \#\hat G=m^3.$$}
\end{thm}

\begin{expl} {\rm Now we consider two-parameter
deformations of Calabi--Yau orbifolds with $h^{1,1}(X)=2$.
Note that for these examples listed in TABLE 2, the groups $G_Q$
are more complicated than the cases $h^{1,1}(X)=1$ discussed
above, but in any event, we have the `duality'' relation: $\# G_Q\times \#\hat G=m^3$.} 

\begin{table}
\begin{center}
\begin{tabular}{|c|c|c|c|} \hline
$m$ & $Q$ & $\hat G$ & $\mbox{generator}$ \\ \hline  
$8$ & $(1,1,2,2,2)$ & $(\Z/8\Z)^2$ & $(1,7,0,0,0),(7,1,0,0,0)$  \\ \hline
$12$ & $(1,1,2,2,6)$ & $(\Z/2\Z)\times(\Z/6\Z)^2$  & $(6,6,0,0,0),
(4,0,2,2,0),(0,4,2,2,0)$ \\ \hline
$12$ & $(1,2,2,3,4)$ & $(\Z/6\Z)^2$ & $(4,2,2,0,0), (0,2,2,0,4)$ \\ \hline
$14$ & $(1,2,2,2,7)$ & $(\Z/7\Z)^2$ & $(2,2,2,2,4),(4,4,4,4,2)$\\ \hline
$18$ & $(1,1,1,6,9)$ & $(\Z/3\Z)\times(\Z/6\Z)^2$ & $(6,6,6,0,0),
(3,3,3,0,9),(3,3,3,3,6)$ \\ \hline
\end{tabular}
\vskip 0.2cm
\caption{Mirror Calabi--Yau threefolds II}
\end{center}
\end{table}

{\rm For any of these Calabi--Yau orbifolds $X$, its mirror Calabi--Yau threefold
is defined by a hypersurface ${\mathcal{W}_Q}=0$ in the weighted projective
space ${\BP}^4(Q)/G$ where  ${\mathcal{W}_Q}$ is of the form:
$${\mathcal{W}}_Q=
Y_1^{m_1}+\cdots + Y_5^{m_5}-m\psi Y_1Y_2\cdots Y_5-
k\phi\begin{cases} Y_1^4Y_2^4\quad & \mbox{for $m=8$} \\
                                   Y_1^6Y_2^6\quad & \mbox{for $m=12$} \\
                                   Y_1^6Y_4^2\quad & \mbox{for $m=12$} \\
                                   Y_1^7Y_5  \quad & \mbox{for $m=14$} \\
                                   Y_1^6Y_2^6Y_3^6 \quad & \mbox {for $m=18$}
                       \end{cases}$$
where $\psi$ and $\phi$ are parameters and $k$ is a positive integer such 
that $k|m$. The mirror of the Calabi--Yau threefolds
with $m=8$ and $m=12,\,Q=(1,1,2,2,6)$ with two parameters was studied 
in Candelas et al. [COFKM93].}
\end{expl}

\begin{rem}
{\rm Theorem 4.2 is valid even when one passes to a pair of finite
groups $(H, \hat H)$ where $H$ is a subgroup of $(\Z/m\Z)^3$ and $\hat H$
is a dual group of $H$ in the sense that $\# H\times \# \hat H=m^3$.
See Klemm and Theisen [KT93] for examples of one-modulus Calabi--Yau threefolds
of this type. See Problem 10.3.} 

\end{rem} 

\begin{rem} {\rm One-(or multi)-parameter deformations of 
Calabi--Yau orbifolds constructed in Theorem 3.1 are no longer 
dominated by product varieties. (Confer Schoen [Sch96].)  
Consequently, these deformations are no longer of CM type varieties.  }
\end{rem}

\section{Fermat motives and mirror maps}\label{Fermatmotivesandmirrormaps}

In this section, we first construct Fermat motives arising from  
the Calabi--Yau threefolds constructed in Theorem 3.1
at special points of deformations, i.e., at the Fermat points (where all
the deformation parameters are set to zero).
We will interpret the mirror symmetry phenomenon for our Calabi--Yau
orbifolds and their mirror partners in terms of Fermat motives. In particular,
we will determine Fermat motives which are invariant under the mirror
map for each mirror pair of Calabi--Yau orbifolds.  

Now we will construct Fermat motives for the Calabi--Yau orbifolds $X$ 
constructed from Fermat hypersurface of degree $m$ and weight $Q=(q_1,q_2,q_3,q_4,q_5)$.
Recall from Section 2 that associated to the Fermat hypersurface threefold
$$V\,:\, Z_1^m+Z_2^m+Z_3^m+Z_4^m+Z_5^m=0\subset\BP^4,$$ 
we have a group
$$\frakG:=\{\bg=(g_1,g_2,g_3,g_4,g_5)\,|\, g_i\in\mu_m\,\}/\{(g,g,\cdots,g)\}\subset
\mbox{Aut}(V)$$ and its dual group 
$$\hat{\frakG}=\{\ba=(a_1,a_2,a_3,a_4,a_5)\,|\, a_i\in(\Z/m\Z),\,\sum_{i=1}^5 a_i\equiv 0\pmod m\,\}.$$
We consider its subset
$$\frakA=\{\ba\in\hat{\frakG}\,|\, a_i\neq 0 \in (\Z/m\Z)\,\forall\, i\,\}.$$
There is a pairing between $\frakG$ and $\frakA$ given by
$$(\bg,\ba)\mapsto \prod_{i=1}^5 g_i^{a_i}\in \bL=\Q(\zeta_m).$$  

Now we pass onto weighted Fermat Calabi--Yau threefolds. We ought to bring in weights
to our discussion.

\begin{thm} {\rm (Yui [Y05])} {\sl Let $Y$ be the singular Calabi--Yau orbifold corresponding to an admissible
pair $<m,Q>$ with $=(q_1,q_2,q_3,q_4,q_5)$ in Theorem 3.1. Define
$$\frakA(Q):=\{\ba=(a_1,a_2,a_3,a_4,a_5)\in\frakA\,|\, a_i\in(q_i\Z/m\Z)\,\forall\,i\,\}.$$
For each character $\ba\in\frakA(Q)$, let $A=[\ba]$ be the $(\Z/m\Z)^{\times}$-orbit of $\ba$,
i.e., $A=\{t\ba\,|\, t\in(\Z/m\Z)^{\times}\,\}$.  Let $p_{\ba}$ and $p_A$
be the projectors defined as in Definition 2.1.  Then $(Y, p_A)=\cM_A$ is the Fermat motive corresponding
to $A$. Furthermore, we have the motivic decomposition
$$(Y,{\bf 1})=\oplus_{A\in O(\frakA(Q))} \cM_A$$
where $A$ runs over the set of $(\Z/m\Z)^{\times}$-orbits in $\frakA(Q)$.

Then the following assertions hold.

(a) We have
$$h^{3,0}(Y)=1=h^{3,0}(\cM_Q)$$
$$h^{2,1}(Y)=\sum_{A\in O(\frakA(Q))} h^{2,1}(\cM_A)$$
and
$$B_3(Y)=\sum_{A\in O(\frakA(Q))} B_3(\cM_A).$$

(b) Let $\pi : X\to Y$ be the crepant resolution of $Y$, and let $\hat{\Sigma}(Y):=\pi^{-1}(\Sigma(Y))$
be the pull-back of the singular locus $\Sigma(Y)$. Let $C\in\Sigma(Y)$ be an
irreducible curve of genus $g(C)$.  Then
$$h^{2,1}(\hat{\Sigma}(Y))=\sum_C m_Cg(C)$$
where the sum runs over all distinct irreducible curves in $\Sigma(Y)$
with multiplicity $m_C$.

(c) For the crepant resolution $X$ of $Y$, we have
$$h^{3,0}(X)=1=h^{3,0}(\cM_A)$$
$$h^{2,1}(X)=h^{2,1}(Y)+h^{2,1}(\hat{\Sigma}(Y)$$
and
$$B_3(X)=B_3(Y)+2\sum_C m_Cg(C).$$

Furthermore,
$$B_2(X)=h^{1,1}(X)=1+\sum_{P}m_P$$
where $P$ is a singular point in  $\Sigma(Y)$ with multiplicity $m_P$.}
\end{thm}  

\begin{expl} 
{\rm Let $m=5$ and $Q=(1,1,1,1,1)$. This is the quintic Calabi--Yau threefold $X$,
which is smooth.  The Fermat motives are tabulated in TABLE 3.

\begin{table}
\begin{center}
\begin{tabular}{|c|c|c|c|c|c|}\hline \hline
A  & $\mbox{dim}$ & $\mbox{mult}$ & $h^{3,0}$ & $h^{2,1}$ & $\sum B_3(\cM_A)$ \\ \hline
$[1,1,1,1,1]$ & $4$    & $1$           & $1$      & $1$       & $4$ \\ \hline
$[1,1,1,3,4]$ & $4$    & $20$          & $0$      & $2$       & $80$ \\ \hline
$[1,1,2,2,4]$ & $4$    & $30$          & $0$      & $2$       & $120$ \\ \hline
\end{tabular}
\vskip 0.2cm
\caption{Fermat motives for the quintic Calabi--Yau threefold}
\end{center}
\end{table}
\smallskip

Here each motive has dimension $\varphi(5)=4$ and multiplicity is given by the number
of permutations of the components.  Furthermore, $Y$ is smooth, so $Y=X$. Summing over 
the motivic Hodge and Betti numbers, we get the global Hodge and Betti numbers:
$$h^{3,0}(X)=1,\,\, h^{2,1}(X)=1+2\times 20+2\times 30=101,\,\, h^{1,1}(X)=1.$$
The third Betti number and the Euler characteristic are given by
$$B_3(X)=4+80+120=204,\,\,E(X)=2(h^{1,1}(X)-h^{2,1}(X))=2(1-101)=-200.$$} 
\end{expl}

\begin{expl}
{\rm Let $m=8$ and $Q=(1,1,2,2,2)$. This is an octic Calabi--Yau threefold $X$. 
The Fermat motives for the singular Calabi--Yau threefold $Y$ are tabulated in TABLE 4.

\begin{table}
\begin{center}
\begin{tabular}{|c|c|c|c|c|c|} \hline \hline
A & $\mbox{dim}$ & $\mbox{mult}$ & $h^{3,0}$ & $h^{2,1}$ & $\sum B_3(\cM_A)$ \\ \hline
$[1,1,2,2,2]$ & $4$ & $1$ & $1$ & $1$ & $4$ \\ \hline
$[7,3,2,2,2]$ & $4$ & $1$ & $0$ & $2$ & $4$ \\ \hline
$[6,4,2,2,2]$ & $2$ & $2$ & $0$ & $2$ & $4$ \\ \hline
$[7,1,2,2,4]$ & $4$ & $6$ & $0$ & $12$ & $24$ \\ \hline
$[6,2,2,2,4]$ & $2$ & $6$ & $0$ & $6$ & $12$ \\ \hline
$[4,4,2,2,4]$ & $2$ & $3$ & $0$ & $3$ & $6$ \\ \hline
$[5,1,2,2,6]$ & $4$ & $3$ & $0$ & $6$ & $12$ \\ \hline
$[4,2,2,2,6]$ & $2$ & $6$ & $0$ & $6$ & $12$ \\ \hline
$[3,3,2,2,6]$ & $4$ & $3$ & $0$ & $6$ & $12$ \\ \hline
$[5,1,2,4,4]$ & $4$ & $3$ & $0$ & $6$ & $12$ \\ \hline
$[4,2,2,4,4]$ & $2$ & $6$ & $0$ & $6$ & $12$ \\ \hline
$[3,3,2,4,4]$ & $4$ & $3$ & $0$ & $6$ & $12$ \\ \hline
$[3,1,2,4,6]$ & $4$ & $6$ & $0$ & $12$ & $24$ \\ \hline
$[2,2,2,4,6]$ & $2$ & $6$ & $0$ & $6$ & $12$ \\ \hline
$[3,1,4,4,4]$ & $4$ & $1$ & $0$ & $2$ & $4$ \\ \hline
$[2,2,4,4,4]$ & $2$ & $1$ & $0$ & $1$ & $2$ \\ \hline
\end{tabular}
\vskip 0.2cm
\caption{Fermat motives for the octic Calabi--Yau threefold}
\end{center}
\end{table}
\smallskip

Note that the dimension of each motive is either $4$ or $2$ and the latter occurs
when $\mbox{gcd}(\ba,8)\neq 1$.  

Summing up the motivic Hodge and Betti numbers, we obtain
$$h^{3,0}(Y)=1,\, h^{2,1}(Y)=83\quad\mbox{and}\quad
B_3(Y)=\sum_{A} B_3(\cM_A)=168.$$
There is one curve singularity $C: Y_3^4+Y_4^4+Y_5^4=0$ which has genus $g(C)=3$ with
multiplicity $m_C=1$. Passing onto the crepant resolution $X$ of $Y$, we then get
$$h^{3,0}(X)=1,\quad h^{2,1}(X)=3+83=86,\quad h^{1,1}(X)=2$$
and
$$B_3(X)=168+2\times 3=174.$$
Finally the Euler characteristic is
$$E(X)=2(2-86)=2(-84)=-168.$$}
\end{expl}

\begin{rem} {\rm This gives an alternative method for computing the global Betti numbers and Euler
characteristic for the Calabi--Yau orbifolds $X$ corresponding to admissible
pairs $<m,Q>$ from the formula of Vafa described in Theorem 3.1 (c).} 
\end{rem}

Now we pass onto mirror partners of Calabi--Yau orbifolds. 
We need to bring in weights into the duality.   

\begin{defn}
{\rm  Let $X$ be a Calabi--Yau orbifold of degree $m$ and weight $Q=(q_1,q_2,q_3,q_4,q_5)$.
Let $\zeta_m$ denote a primitive $m$-th root of unity.
Let $m_i=m/q_i$ and let $\mu_{m_i}$ denote the group (scheme) of
$m_i$-th roots of unity, for each $i, \,0\leq i\leq 5$.
Let
$$\frakA(Q)=\{\,\ba=(a_1,a_2,a_3,a_4,a_5)\in\frak A\,|\, a_i\in(q_i\Z/m\Z)\,\}.$$ 
Define 
$$\hat G=\{\,\bg=(g_1,g_2,g_3,g_4,g_5)\,|\, g_i\in\mu_{m_i},\, 
\prod_{i=1}^5 g_i=1\,\}/\{(g,g,\cdots,g)\,\}.$$ 
Then $\hat G$ is the group of discrete symmetry on $X$ leaving the holomorphic $3$-form
$\Omega$ invariant, and $\hat G$ acts on the set $\frakA(Q)$ by
$$\hat G\times \frakA(Q)\to \Q(\zeta_m)\,:\,
(\bg,\ba)\mapsto \prod_{i=1}^5 g_i^{a_i}=g_1^{a_1}g_2^{a_2}\cdots g_5^{a_5}.$$
\smallskip

(1) Let $\ba\in\frakA(Q)$, and $\bg\in\hat G$ be a generator of $\hat G$. 
We say that $\ba$ is {\it invariant under $\bg$} if 
$$\prod_{i=1}^5 g_i^{a_i}=g_1^{a_1}g_2^{a_2}g_3^{a_3}g_4^{a_4}g_5^{a_5}=1.$$
\smallskip

(2) Let $A=[\ba]$ be the $(\Z/m\Z)^{\times}$-orbit of $\ba$ and 
let $\cM_A$ be the corresponding Fermat motive.
We say that $\cM_A$ is {\it invariant under $\hat G$} if any $\ba\in A$ is
invariant under every generator $\bg$ of $\hat G$.}
\end{defn}

\begin{expl}
{\rm We consider Calabi--Yau orbifolds $X$ and their mirror partners with
K\"ahler modulus $1$, i.e., $h^{1,1}(X)=1$. There are altogether
four such Calabi--Yau orbifolds among the $147$ cases.
All these four Calabi--Yau orbifolds, no singularities of dimension $1$.
\smallskip

(a) Let $m=5$ and $Q=(1,1,1,1,1)$. We know that
$h^{2,1}=101$, $B_3=204$, and that 
$\hat G=(\Z/5\Z)^3$ generated by 
$$\bg=(1,0,0,4,0),(1,0,4,0,0),\quad\mbox{and}\quad (1,4,0,0,0).$$
Let $\ba=(1,1,1,1,1)$.  Then for $\bg=(1,0,0,4,0)$, we have 
$$\prod_{i=1}^5 g_i^{a_i}=(\zeta_5)^1\cdot 1^1\cdot 1^1\cdot
(\zeta_5^4)^1\cdot 1^1=\zeta_5^5=1.$$ 
and similarly for the other two generators. 
Now consider the Fermat motive corresponding to
$$[Q]=[(1,1,1,1,1)]=\{(1,1,1,1,1),(2,2,2,2,2),(4,4,4,4,4),(3,3,3,3,3)\}.$$
We can compute that all elements in $[Q]$ are invariant under $\hat G$.
Hence the motive $\cM_{Q}$ is invariant under $\hat G$.

On the other hand, let $\ba=(1,1,1,3,4)$. Then for $\bg=(1,0,0,4,0)$, we have
$$\prod_{i=1}^5g_i^{a_i}=(\zeta_5)^1\cdot 1^1\cdot 1^1\cdot (\zeta_5^4)^3\cdot 1^4 
=(\zeta_5)^{13}\neq 1,$$ 
and similarly for the other two generators.
Also for $\ba=(1,1,2,2,4)$ and $\bg=(1,0,0,4,0)$, we have 
$$\prod_{i=5}^5g_i^{a_i}
=(\zeta_5)^1\cdot 1^1\cdot 1^2\cdot (\zeta_5^4)^2\cdot 1^4=(\zeta_5)^9\neq 1,$$
and similarly for the other two generators.
Consequently, the Fermat motives $\cM_A$ for $A=[(1,1,1,3,4)],\,[(1,1,2,2,4)]$
are not invariant under $\hat G$.
\smallskip

(b) Let $m=6$ and $Q=(1,1,1,1,2)$. We know that 
$h^{2,1}=103$, $B_3=208$ and that $\hat G=(\Z/3\Z)\times (\Z/6\Z)^2$
with generators 
$$\bg=(0,2,2,2,0), (5,0,0,1,0), (0,5,0,1,0).$$ 
The Fermat motives for this Calabi--Yau orbifold are:
$$(1,1,1,1,2),(1,1,1,5,4),(1,1,2,4,4),(1,1,3,3,4),(1,1,3,5,2),$$
$$(1,1,4,4,2),(1,2,2,3,4),(1,2,2,5,2),(1,2,3,4,2),(1,3,3,3,2),$$
$$(2,2,2,2,4),(2,2,2,4,2),(2,2,3,3,2)$$
and they are all of dimension $2$.
 
Among these, the Fermat motives corresponding to $(1,1,1,1,2)$ and $(2,2,2,2,4)$
are invariant under $\hat G$. For instance, for $\ba=(2,2,2,4,2)$ and $\bg=(5,0,0,1,0)$,
we have
$$\prod_{i=1}^5 g_i^{a_i}=(\zeta_6^5)^2\cdot 1^2\cdot 1^2\cdot (\zeta_6)^2\cdot 1^2
=(\zeta_6)^{12}=1,$$
and similarly for the other generators.
On the other hand, for $\ba=(1,1,3,3,4)$ and $\bg=(5,0,0,1,0)$, we compute
$$\prod_{i=1}^5 g_i^{a_i}=(\zeta_6)^8\neq 1,$$ and also for the other
two generators. 
Therefore, the two Fermat motives, 
$\cM_A$ with $A=[Q]=[(1,1,1,1,2)],\,[(2,2,2,2,4)]$ are invariant under
$\hat G$.
\smallskip

(c) Let $m=8$ and $Q=(1,1,1,1,4)$. Then $h^{2,2,1}=149$, $B_3=300$, and that
$\hat G=(\Z/2\Z)\times(\Z/8\Z)^2$ with generators
$\bg=(7,0,0,1,0),(0,7,0,1,0),(0,0,7,1,0)$. After calculating the action
of $\hat G$ on each Fermat motive, we see that only the motive
corresponding to weight $Q$ of dimension $4$ is invariant. 
\smallskip

(d) Let $m=10$ and $Q=(1,1,1,2,5)$. We know that $h^{1,1}=1=B_2$, $h^{2,1}=145$
and $B_3=284$, and that $\hat G=(\Z/10\Z)^2$ with generators
$\bg=(9,0,1,0,0),(0,9,1,0,0)$. 
We compute the products $\prod_{i=1}^5 g_i^{a_i}$ for all motives and generators
and the outcome is that only the motive $\cM_{Q}$ corresponding to the
weight $Q$ of dimension $4$ is invariant under $\hat G$.

In each case of (a),(c) and (d), the Fermat motive $\cM_Q$ corresponding to the weight
$Q$ is the only motive left invariant under the mirror map. It is of dimension $4$, and we have
$B_3(\hat X)=4=2(1+h^{1,1}(X))$, while, in the case (b), there are two Fermat motives
$\cM_Q$ and $\cM_{[2,2,2,2,4]}$ which are invariant under the mirror map, and
we have
$B_3(\hat X)=2+2=2(1+h^{1,1}(X)).$}
\end{expl}

\begin{expl}
{\rm We consider examples of Calabi--Yau orbifolds $X$ and their
mirror partners with K\"aher modulus two, i.e., $h^{1,1}(X)=2$.
There are six such examples among our $147$ Calabi--Yau orbifolds.
Some of these six Calabi--Yau orbifolds have curve singularities.
\smallskip

(a) Let $m=8$ and $Q=(1,1,2,2,2)$. We know that $h^{2,1}=86$
and $B_3=168+(6)$ where $6$ comes from the singular locus, 
and that $G_Q=(\Z/2\Z)^3$ and $\hat G=(\Z/8\Z)^2$ with
generators $\bg=(1,7,0,0,0),(7,1,0,0,0)$. The Fermat motives constituting
singular part of this Calabi--Yau orbifold correspond to:
$$(1,1,2,2,2),(1,1,2,6,6),(1,1,4,4,6),(1,3,2,4,6),(1,3,2,2,6),(1,3,4,4,4),$$
$$(1,5,2,2,6),(1,5,2,4,4),(1,5,6,6,6),(1,7,2,2,4),$$
of dimension $4$, and
$$(2,2,2,4,6),(2,2,4,4,4),(2,4,2,2,6),(2,4,2,4,4),(2,6,2,2,4),(4,4,2,2,4)$$
of dimension $2$.

We compute the product $\prod_{i=1}^5 g_i^{a_i}$ for all these motives
and all generators $\bg$.  For the weight $Q=(1,1,2,2,2)$ and $\bg=(1,7,0,0,0)$,
we have
$$\prod_{i=1}^5 g_i^{a_i}=(\zeta_8)^1\cdot(\zeta_8^7)^1\cdot 1^2\cdot 1^2\cdot 1^2=1$$
and similarly for $\bg=(7,1,0,0,0)$. Therefore the motive $\cM_{Q}$ is
invariant under $\hat G$. The motive $\cM_{Q}$ has dimension $4$. 
Also for $(2,2,4,4,4)$ and $\bg=(1,7,0,0,0)$, we have
$$\prod_{i=1}^5 g_i^{a_i}=(\zeta_8)^2\cdot(\zeta_8^7)^2\cdot 1^4\cdot 1^4\cdot 1^4=1$$
and similarly for $\bg=(7,1,0,0,0)$. The motive $\cM_{[2,2,4,4,4]}$ has dimension $2$.
All the remaining Fermat motives are not invariant under $\hat G$. 


\smallskip

(b) For $m=12$ and $Q=(1,1,2,2,6)$, we know that $h^{2,1}=128$ and $B_3=258$. 
The weight motive $\cM_{[1,1,2,2,6]}$ has dimension $4$. There is another motive of dimension $2$
which is invariant under the mirror operation. 
\smallskip

(c) For $m=12$, $Q=(1,2,2,3,4)$, we know that $h^{2,1}=74$, $B_3=150$. 
The weight motive $\cM_{[1,2,2,3,4]}$ has dimension $4$ and it is
the only motive which is invariant under the mirror operation. 

(d) Let $m=14$ and $Q=(1,2,2,2,7)$. Then $h^{2,1}=122$, $B_3=216+(30)$ where $30$
comes from the singularity, and $\hat G=(\Z/7\Z)^2$ with generators
$\bg=(2,2,2,2,4)$ and $(4,4,4,4,2)$.  We calculate that only the motive
$\cM_{Q}$ corresponding to weight $Q$ is invariant under $\hat G$,
and it has dimension $6$.
\smallskip

(e) Let $m=18$ and $Q=(1,1,1,6,9)$. Then $h^{2,1}=272$ and $B_3=546$, and 
$\hat G=(\Z/3\Z)\times(\Z/6\Z)^2$ with generators
$\bg=(6,6,6,0,0), (3,3,3,0,9), (3,3,3,3,6)$. We calculate that there is
only one motive that is left invariant under $\hat G$, namely,
the Fermat motive $\cM_{Q}$ corresponding to the weight $Q$, and its
dimension is $6$.} 
\end{expl}

\begin{thm} {\sl Let $(X, \hat X)$ be a mirror pair of Calabi--Yau orbifold corresponding to
an admissible pair $<m, Q>$ with $Q=(q_1,q_2,q_3,q_4,q_5)$. Let
$\bg=(g_1,g_2,g_3,g_4,g_5)$ be any generator of $\hat G$. Let
$\cM_A=[a_1,a_2,a_3,a_4,a_5]$ be a Fermat motive. 

Suppose that
$$\prod_{i=1}^5 g_i^{a_i}=1\,\,\forall\,\,\ba\in A$$
then $\cM_A$ is invariant under the mirror map.

In particular, if $h^{1,1}(X)=1$, then the Fermat motive $\cM_Q$ corresponding to
the weight $Q$ is the only Fermat motive left invariant under $\hat G$.

If $h^{1,1}(X)>1$, there is an algorithm to determine the other Fermat motives
which are invariant under the mirror map.}
\end{thm}

{\rm This theorem will be proved in Section 9 below. We ought to understand
possible relations between motives and monomials and integral points in reflexive
polytopes in toric geometric setting \`a la Batyrev [Ba94].}

\begin{rem} {\rm When $h^{1,1}(X)=1$, Goto [G06] has recently showed that
the formal group arising from $X$ is invariant under the mirror map. In particular,
the height $h$ of the formal group at bounded above by $2$.}
\end{rem}

\section{Batyrev's mirror symmetry}\label{Batyrevmirrorsymmetry}

Batyrev [Ba94] gives a combinatorial construction of mirror pairs of
Calabi--Yau hypersurfaces in the toric geometric setting. Here is a brief summary of
Batyrev's construction and the main ingredients in his theory.  A nice reference
on this topic might be Cox and Katz [CK99]. 

Let $\Delta\subset\R^n$ be an $n$-dimensional polytope, and let $\Delta^*=\mbox{Hom}(\Delta,\Z)$
be the dual polytope. Denote by $<*,*>$ the nondegenerate pairing between the $n$-dimensional
$\R$-vector spaces $\Delta_{\R}$ and $\Delta^*_{\R}$.  Then 
$$\Delta^*:=\{{\mathbf y}=(y_1,y_2,\cdots, y_n)\,|\, <{\mathbf y},{\mathbf x}>:=\sum_{i=1}^n y_ix_i\geq -1\,\forall\, 
{\mathbf x}=(x_1, x_2,\cdots, x_n)\in\Delta\,\}.$$
The polytope $\Delta$ is said to be {\it reflexive} if it has the following properties:
\smallskip

(1) $\Delta$ is convex integral polytope, i.e., all vertices of $\Delta$ are integral,

(2) $\Delta$ contains the origin $v_0=(0,0,\cdots,0)$ as an interior point, and

(3) each codimension one face is of the form 
$$\{\,{\mathbf x}\in{\Delta}_{\R}\,|\,< {\mathbf y}, {\mathbf x} >=-1\,\,\mbox{for some}
\,\,{\mathbf y}\in\Delta^*_{\R}\,\}.$$
\smallskip

If $\Delta$ is reflexive, then its dual $\Delta^*$ is again reflexive, and $(\Delta^*)^*=\Delta$.
(Indeed, this fact is the basis for Batyrev's mirror symmetry.)

To a reflexive polytope $\Delta$, one associates a complete rational fan $\Sigma(\Delta)$ as
follows: For every $l$-dimensional face $\Theta_{l}\subset\Delta$, define an
$n$-dimensional cone $\sigma(\Theta_{l})$ by
$$\sigma(\Theta_{l}):=\{\lambda(p^{\prime}-p)\,|\, \lambda\in\R_{\geq 0}, p\in\Delta,\, p^{\prime}\in\Theta_{\ell}\,\}.$$
That is, $\sigma(\Theta_{l})$ consists of all vectors $\lambda(p-p^{\prime})$ where
$\lambda\in\R_{\geq 0},\, p\in\Delta,\,p^{\prime}\in\Theta_l$. Then the fan $\Sigma(\Delta)$ is given as 
the collection of all $(n-l)$-dimensional dual cones
$\sigma^*(\Theta_{l})$ for $l=0,1,\cdots, n$ for all faces of $\Delta$, and the complete fan defines the toric variety
$\BP_{\Delta}$. (That is, $\BP_{\Sigma(\Delta)}$ is defined by a compactification of the algebraic torus $(\C^*)^n$ using
combinatorial data encoded in the fan $\Sigma(\Delta)$.)

\begin{lem} {\rm (Batyrev [Ba94])} {\sl Denote by $v_i\,\, (i=0,1,\cdots, s)$ the integral points in $\Delta$. 
Consider the affine space $\C^{s+1}$ with
smooth coordinates ${\mathbf c}=(c_0,c_1,\cdots, c_s)$.  Let $Z_f$ denote
the zero locus of the Laurent polynomial
$$f_{\Delta}({\mathbf c}, {\mathbf X})=\sum_{i=0}^s c_i {\mathbf X}^{\mathbf{\mu}}\in\C[X_1^{\pm 1}, X_2^{\pm 1},\cdots,
X_n^{\pm 1}]$$
where ${\mathbf X}^{\mathbf{\mu}}=X_1^{\mu_1}X_2^{\mu_2}\cdots X_n^{\mu_n}$ in the algebraic torus $(\C^*)^n\subset
\BP_{\Delta}$. Let $\bar{Z}_f$ be the closure of $Z_f$ in $\BP_{\Delta}$.  
The $\Delta$-regularity conditions for hypersurfaces imply that the
singularities of hypersurfaces are induced only by singularities of the ambient toric variety
$\BP_{\Delta}$.  As a consequence, one obtains a simultaneous resolution of all
members of the family of Calabi--Yau hypersurfaces, and hence 
one has a different family of smooth Calabi--Yau hypersurfaces $X_{\Delta}$ associated to $\Delta$.

Similarly, starting with the dual reflexive polytope $\Delta^*$ and carrying out the same construction as above, one
obtains a family of smooth Calabi--Yau hypersurfaces $X_{\Delta^*}$ associated to $\Delta^*$. } 
\end{lem}

We thus obtain a pair of different families of smooth Calabi--Yau hypersurfaces. The remarkable theorem of Batyrev is formulated
as follows.

\begin{thm} {\rm (Batyrev)} {\sl If $n\leq 4$, there is a crepant resolution $X_{\Delta}$ (resp.
$X_{\Delta^*)}$ of singularities for the hypersurface $\bar{Z}_{f_{\Delta}}$ (resp. $\bar{Z}_{f_{\Delta^*}}$). 

When $n=4$, the resolutions $X_{\Delta}$ and $X_{\Delta^*}$ are mirror symmetric in the sense of
Topological Mirror Symmetry Conjecture 1.1, that is,
$$h^{1,1}(X_{\Delta})=h^{2,1}(X_{\Delta^*}),\quad h^{2,1}(X_{\Delta})=h^{1,1}(X_{\Delta^*}).$$

The Hodge numbers are given combinatorially by
$$h^{1,1}(X_{\Delta})=h^{2,1}(X_{\Delta^*})=l(\Delta^*)-(4+1)-\sum_{\mbox{codim}\Theta^*=1} l^{\prime}(\Theta^*)
+\sum_{\mbox{codim}\Theta^*=2}l^{\prime}(\Theta^*)l^{\prime}(\Theta),$$ 
$$h^{2,1}(X_{\Delta})=h^{1,1}(X_{\Delta^*})=l(\Delta)-(4+1)-\sum_{\mbox{codim}\Theta=1}l^{\prime}(\Theta)
+\sum_{\mbox{codim}\Theta=2}l^{\prime}(\Theta)l^{\prime}(\Theta^*).$$
Here $\Theta, \Theta^*$ are faces of $\Delta$ and $\Delta^*$, and $(\Theta,\Theta^*)$ is a dual pair. 
$l(\Delta)$ denotes the number of integral points in $\Delta$, and
$l(\Theta)$ (resp. $l^{\prime}(\Theta)$) the number of integral points (resp. the number of integral 
points in the interior) of the face $\Theta$.}
\end{thm}

\begin{rem} {\rm The formulae for Hodge numbers give yet another method of
computing Hodge numbers (and hence Betti numbers) and the Euler characteristic
of our Calabi--Yau orbifolds.  Via Batyrev's combinational approach, there
are altogether 473 800 776 families of Calabi--Yau threefolds corresponding to
$4$-dimensional reflexive polytopes. See Kreuzer and Skarke [KS ].} 
\end{rem}

Our mirror pairs of Calabi--Yau hypersurfaces obtained in Section 4 by orbifolding construction can
be reformulated in the framework of Batyrev's mirror construction. Confer Batyrev [Ba94]. 

\begin{prop}{\sl Let $<m,Q>$ be an admissible pair where $Q=(q_1,q_2,q_3,q_4,q_5)$ satisfies
the Calabi--Yau condition: $m=q+_1+q_2+\cdots+q_5$, or equivalently $\sum_{i=1}^5\frac{q_i}{m}=1$ given in
Theorem 3.1. Let
$$\Delta(Q)=\{(x_1,x_2,x_3,x_4,x_5)\in\Q^5\,|\,\sum_{i=1}^5 q_ix_i=0,\, x_i\geq -1\,\, (1\leq i\leq 5)\,\}.$$
Then $\Delta(Q)$ is a reflexive polytope and the toric variety $\BP_{\Delta^*(Q)}$ is isomorphic
to the weighted projective $4$-space $\BP^4(Q)$, and the Calabi--Yau threefold 
$X_{\Delta(Q)}$ is isomorphic to some Calabi--Yau orbifolds in Theorem 4.1. 

More precisely, we have the following assertions.

(a) Put $m_i=m/q_i$ for each $i, \,1\leq i\leq 5$. Then the quasi-homogeneous equation
$$Y_1^{m_1}+Y_2^{m_2}+Y_3^{m_3}+Y_4^{m_4}+Y_5^{m_5}=0$$
defines a $\Delta(Q)$-regular Calabi--Yau hypersurface of Fermat type in the
weighted projective $4$-space $\BP^4(Q)$. Then the family of Calabi--Yau hypersurfaces consists
of quotients of deformations of this hypersurface by the fundamental group $\pi_1(\Delta(Q))$.
Here $\pi_1(\Delta(Q))$ is isomorphic to the kernel of the surjective homomorphism
$$(\mu_{m_1}\times\mu_{m_2}\times\mu_{m_3}\times\mu_{m_4}\times\mu_{m_5})/\mu_n\to\mu_m$$
and it has order $m^3/q_1q_2\cdots q_5$.

(b) For the dual reflexive polytope $\Delta^*(Q)$, the fundamental group
$\pi_1(\Delta^*(Q))$  is isomorphic to the kernel of the surjective
homomorphism
$$(\mu_{m_1}\times\mu_{m_2}\times \cdots\times \mu_{m_5})/\mu_m\to \mu_m\quad:\quad
\bar{\gamma}_Q(g_1^{a_1}\cdots g_5^{a_5})=g^{q_1a_1+\cdots q_5a_5}$$
and it has order $q_1q_2q_3q_4q_5=\# G_Q$.  

(c) In the toric geometric setting, mirror pairs of Calabi--Yau hypersurfaces
is described by a pair of reflexive polytopes, that is,
$\pi_1({\Delta(Q)})$ and $\pi_1({\Delta^*(Q)})$ are dual finite abelian groups
in the sense that
$$\#\pi_1(\Delta(Q))\times \#\pi_1(\Delta^*(Q))=m^3$$
coinciding with the mirror duality in the orbifolding construction of Theorem 4.1.

(d) The mirror map for the families of Calabi--Yau hypersurfaces is the following
Galois correspondence: If $\mathcal F(\Delta(Q))$ and $\mathcal F(\Delta^*(Q))$ are
quotients respectively by $\pi_1(\Delta(Q))$ and $\pi_1(\Delta^*(Q))$ of some subfamilies of
deformations of Calabi--Yau hypersurfaces obtained by the construction in Lemma 6.1.}
\end{prop}

Reflexive polytopes for Calabi--Yau
hypersurfaces of Fermat type are discussed in a number of
literatures listing integral vertices, e.g., Batyrev [Ba94],
Hosono-Lian and Yui [HLY95]. Here are some examples. 

\begin{expl} {\rm In our examples of one-parameter
Calabi--Yau orbifolds with $h^{1,1}=1$ in Example 4.1, we have $q_1=1$
and $q_2\leq q_3\leq q_4\leq q_5$.  The convex polytope
$\Delta(Q)$ is the convex hull of the integral vertices given as follows:
$$(0,0,0,0),(-1,-1,-1,-1),((\frac{m}{q_2})-1,-1,-1,-1),(-1,(\frac{m}{q_3}-1),-1,-1)$$
$$(-1,-1,(\frac{m}{q_4}-1),-1),(-1,-1,-1,(\frac{m}{q_5}-1)).$$
The dual polytope $\Delta^*(Q)$ consists of the following vectors:
$$(0,0,0,0), (1,0,0,0),(0,1,0,0),(0,0,1,0),(0,0,0,1),(-q_2,-q_3,-q_4,-q_5).$$
}
\end{expl}

\begin{expl} {\rm For the two-parameter families of Calabi--Yau orbifolds
with $h^{1,1}=2$ in Example 4.2, we again have $q_1=1$ and $q_2\leq q_3\leq q_4\leq q_5$.
We list the dual polytopes for these examples in TABLE 5. }

\medskip

\begin{table}
\begin{center}
\begin{tabular}{|c|c|c|} \hline
$m$ & $Q$ & $\Delta^*(Q)$ \\ \hline
$8$ & $(1,1,2,2,2)$ & $(0,0,0,0),(1,0,0,0),(0,1,0,0),(0,0,1,0),$ \\  
$\quad$ & $\quad$ & $(0,0,0,1),(-1,-2,-2,-2),(-1,-1,-1,0)$ \\ \hline
$12$ & $(1,1,2,2,6)$ & $(0,0,0,0),(1,0,0,0),(0,1,0,0),(0,0,1,0),$ \\
$\quad$ & $\quad$ &  $(0,0,0,1),(-1,-2,-2,-6),(-3,-1,-1,0)$ \\ \hline
$12$ & $(1,2,2,3,4)$ & $(0,0,0,1),(1,0,0,0),(0,1,0,0),(0,0,1,0),$ \\
$\quad$ & $\quad$    & $(0,0,0,1),(-2,-2,-3,-4),(-1,-1,-1,-2)$ \\ \hline
$14$ & $(1,2,2,27)$ & $(0,0,0,0),(1,0,0,0),(0,1,0,0),(0,0,1,0),$ \\
$\quad$ & $\quad$ & $(0,0,0,1),(-2,-2,-2,-7),(-1,-1,-1,-3)$ \\ \hline
$18$ & $(1,1,1,6,9)$ & $(0,0,0,0),(1,0,0,0),(0,1,0,0),(0,0,1,0),$ \\
$\quad$ & $\quad$ & $(0,0,0,1),(-1,-1,-6,-9), (0,0,-2,-3)$ \\ \hline
\end{tabular}
\vskip 0.2cm
\caption{Vertices for Dual Polytopes} 
\end{center}
\end{table}
\end{expl}

\begin{prop}
{\sl For our mirror paris of Calabi--Yau hypersufaces constructed in
Theorem 4.1 and Proposition 6.3, the origin is the only integral point
present in both the reflexive polytope $\Delta(Q)$ and its dual reflexive
polytope $\Delta^*(Q)$.}
\end{prop}


In Batyrev's mirror symmetry, the essential information is encoded in the {\it Newton polyhedron} of the
defining equation of the Calabi--Yau hypersurface, which is the convex hull of the monomials appearing in the
hypersurface.   
Suppose that a family of Calabi--Yau hypersurfaces in $\BP^4(Q)$ is the zero locus
of the generic equation 
$$f({\mathbf c},{\mathbf X})=0\subset\BP^4(Q)$$
of degree $m$ and weight $Q=(q_1,q_2,q_3,q_4,q_5)$ with deformation parameters.
The Newton  polyhedron of $f$ always contains
the monomial $x_1\cdots x_5$. Corresponding to this monomial, we have the integral point $(1,1,\cdots, 1)$. Translate it to the
origin. Then the Newton polyhedron of $f$ is given by 
$$\Delta(Q)=\mbox{convex hull of}\,\{(x_1,x_2,x_3,x_4,x_5)\in\Q^5\,|\, 
\sum_{i=1}^5 q_ix_i=0,\, x_i\geq -1 \,(1\leq i\leq 5)\,\}$$
The shifted by $(1,1,\cdots,1)$ convex polytope is the intersection of $\Q^5$ with the affine hyperplane
$$q_1x_1+q_2x_2+\cdots+q_5x_5=q_1+q_2+\cdots+q_5=m.$$
Since $X_{\Delta(Q)}$ is isomorphic to some of our Calabi--Yau orbifolds in
Theorem 4.1, we obtain the following result.

\begin{thm} {\sl 
There is a one-to-one correspondence between integral points in $\Delta(Q)$ 
and monomials of degree $m$ in the graded polynomial ring $\mathcal R=\C[X_1,X_2,X_3,X_4,X_5]$ of 
a Calabi--Yau orbifold $X$ associated to an admissible pair $<m,Q>$.  }
\end{thm}

\begin{rem} {\rm
In fact, this is a concrete realization of the more general 
high-power monomial--divisor mirror map constructed 
by Batyrev [Ba94] and Aspinwall--Greene and Morrison [AGM93] giving the isomorphism between 
the subspaces $H^{1,1}_{toric}(X)$ and $H^{2,1}_{poly}(\hat X)$.  
The space $H^{2,1}_{poly}(\hat X)\subset H^{2,1}(\hat X)$ is isomorphic to 
the first-order polynomial deformations, and can be generated by monomials. Its
dimension $h^{2,1}_{poly}(\hat X)$ is equal to the number of generators of
the Mori cone (which is dual to the K\"ahler cone of the mirror $X$). 
The generators of the Mori cone correspond to divisor classes which generate
the subspace $H^{1,1}_{toric}(X)$. This gives rise to a natural isomorphism
$$H^{1,1}_{toric}(X)\simeq H^{2,1}_{poly}(\hat X)$$
and such an isomorphism can be interpreted as the differential of the expected
mirror map between the moduli spaces.}
\end{rem}  


\section{Monomials and periods}\label{monomialsandperiods}

We will briefly describe the Dwork--Katz--Griffiths reduction method [CK99] 
for determining the Picard--Fuchs equations for hypersurfaces 
in a weighted projective space $\BP^4(Q)$ with weight $Q=(q_1,q_2,q_3,q_4,q_5)$. 
The Picard--Fuchs differential equation is the device which establish the one-to-one 
correspondence between monomials and periods. 

Suppose that a hypersurface $P({\mathbf X}):=P(X_1,\cdots, X_5)=0$ defines
a Calabi--Yau threefold $X$ of degree $m$ and weight $Q=(q_1,q_2,\cdots,q_5)$, and let 
$P_{\mathbf c}({\mathbf  X})$ be a Calabi--Yau threefold $X_{\mathbf c}$ with a deformation 
parameter ${\mathbf c}$ (in practice, these parameters are complex structure moduli.)
Let $\Omega_{\mathbf c}$ denote the unique (up to a scalar multiplication) holomorphic
$3$-form on $X_{\mathbf c}$. Then the period  
$\Pi({\mathbf c})$ can be written as
\begin{equation}
 \Pi({\mathbf c})=\int_{\Gamma_i}\Omega_{\mathbf c}=\int_{\gamma}\int_{\Gamma}\frac{\omega}
{P_{\mathbf c}}
\end{equation}
where
$$\omega=\sum_{i=1}^{5}(-1)^{i}q_{i}X_{i}dX_1\wedge\ldots\wedge \hat{dX_i}\wedge\ldots\wedge 
dX_5$$
is the unique (up to scalar) holomorphic $3$ form on $X_{\mathbf c}$. Here
$\Gamma_i\in H_3(X_{\mathbf c}, \Z)$ is a topological $3$-cycle, 
and $\gamma$ is a small curve around the hypersurface
$P=0$ in the $4$-dimensional embedding space. 
Since $H_3(X,\Z)$ has rank $B_3(X)$, there are in total $B_3(X)$ periods.

\begin{defn} {\rm
Let $\C[X_1,X_2,X_3,X_4,X_5]$ be the weighted polynomial ring with weight
$Q=(q_1,q_2,q_3,q_4,q_5)$.  Let ${\mathbf X}^{\mathbf v}=X_1^{v_1}X_2^{v_2}\cdots X_5^{v_5}$
be a monomial in $\C[X_1,X_2,\cdots, X_5]$.
The {\it degree} of ${\mathbf X}^{\mathbf v}$, denoted by $w({\mathbf v})$ is defined by by
$$w({\mathbf v})=\sum_{i=1}^5 q_iv_i=q_1v_1+q_2v_2+\cdots +q_5v_5.$$}
\end{defn}

We see that there are in total $B_3(X)$ monomials, which
are divided into monomials of degree $w({\mathbf v})=0, m,\, 2m,\, 3m$ with
cardinalities $1, \, h^{2,1},\, h^{1,2},\, 1$, respectively.  

\begin{prop}{\sl 
Let ${\mathbf X}^{\mathbf v}=X_1^{v_1}X_2^{v_2}\cdots X_5^{v_5}$ be a monomial
of degree $w({\mathbf v})=\sum_{i=1}^5 q_iv_i$ in the
weighted polynomial ring $\C[X_1,X_2,X_3,X_4,X_5]$ with weight
$(q_1,q_2,\cdots, q_5)$.  Then there is a one-to-one correspondence
between monomials and periods. That is, to every monomial ${\mathbf X}^{\mathbf v}$,
there corresponds a period
$$\Pi_{{\mathbf v}}=
\int_{\Gamma_{\mathbf v}}\frac{{\mathbf X}^{\mathbf v}}{P_{\mathbf c}^{w({\mathbf v})+1}}$$
where $\Gamma_{\mathbf v}$ is a topological cycle corresponding to ${\mathbf v}$ in
$H_3({\mathbf X}_{\mathbf c},\Z)$.}
\end{prop}

Our examples of mirror pairs of Calabi--Yau threefolds in weighted projective  $4$-spaces are
equipped with large automorphism groups. If $(X, \hat X)$ is a mirror pair of families
of Calabi--Yau threefolds corresponding to $<m,Q>$, then the automorphism group
consists of $\hat G$ and permutations of the weighted projective coordinates. 
We may choose a set of basis elements for the set of monomials under this action.
 
The action of the automorphisms on monomials can be transfered by Proposition 7.1 to the actions of
partial differential operators on periods, and under these actions modulo
the exact ones, periods decompose into the disjoin union of orbits.    

Observe that ${\frac{\partial}{\partial X_i}}\left(\frac{g({\mathbf X})}{P_{\mathbf c}^r}\right)
\omega$ 
is exact if $g({\mathbf X})$ is homogeneous with degree such that the whole expression has degree
zero. This leads to the partial integration rule
$$\frac{h\,\partial_i P}{P_{\mathbf c}^r}=\frac{1}{r-1}\,\frac{\partial_i h}
{P_{\mathbf c}^{r-1}}$$  
with $\partial_i=\frac{\partial}{\partial X_i}$.
One then takes derivatives of the expressions
$$
\Pi_{\mathbf v}=\int_{\Gamma}\frac{{\mathbf X}^{\mathbf v}}{P_{\mathbf c}^{w({\mathbf v})+1}}
$$
with respect to the moduli ${\mathbf c}$ of the defining equation. 

If one produces an expression such that the numerator
in the integrand is not one of the basis elements, one relates it,
using the equations $\partial_iP=\ldots$, to the basis and uses
the relation given above. This leads to a system of first order
differential equations (known as `Gauss--Manin equations') for the
$\Pi_{\mathbf v}$ which can be rewritten as a system of partial
differential equations for the period, which is the Picard--Fuchs
equations: 
$$\partial_{c_k}\Pi({\mathbf c})=M^{(k)}(a)\Pi({\mathbf c}),\;
k=1,\ldots, h^{2,1}.$$ 

Furthermore, the Picard--Fuchs differential equations (or equivalently, the periods) are equipped with the
following operators: 
$${\mathcal D}_i\left(\frac{{\mathbf X}^{\mathbf v}}{P_{\mathbf c}^{w({\mathbf v})+1}}\right)
:=\partial_i\left( X_i\,\frac{{\mathbf X}^{\mathbf v}}{P_{\mathbf c}^{w({\mathbf v})+1}}\right)$$
for $i=1,2,\cdots, 5$. Since
$$\int_{\Gamma} d^5{\mathbf X}{\mathcal D_i}\left(\frac{{\mathbf X}^{\mathbf v}}{P_{\mathbf c}^{w({\mathbf v})+1}}\right)=0,$$
this gives rise to relations (identities) between
the differential form associated to the monimial
${\mathbf X}^{\mathbf v}$ and those associated to ${\mathbf X}^{\mathbf v}$ multiplied by
$\hat G$-invariant monomials.  Hence under the action of ${\mathcal D}_i$ ($i=1,2,\cdots, 5$),
the Picard--Fuchs differential equation (and hence the period) decomposes into
the product of differential equations.
\smallskip

Summing up, the Picard--Fuchs differential equation is determined by those complex structure parameters
for which there is a monomial perturbation in the defining equation, and there are
$h^{2,1}_{poly}$ for these. The number of deformation parameters is equal to the
number of generators of the Mori cone (which is dual to the K\"ahler cone) in the
mirror. The generators of the Mori cone come from relations between the divisors,
in this way we get the secondary fan which describes the complex structure moduli space
in terms of ``large structure coordinates''.

\smallskip 

Based on the theory of Dwork--Katz--Griffiths, these
operators are explicitly defined for the one-parameter family of
the quintic Calabi--Yau threefolds by Candelas et al. [CORV00, CORV03], and for the
two-parameter octic family of Calabi--Yau threefolds by Kadir [Ka04, Ka05]. Also they have
the so-called Picard--Fuchs diagrams to illustrate the decomposition of the Picard--Fuchs
differential equation and the periods. Clearly, such decompositions ought to 
have origin in the motivic decomposition of the manifolds.
We obtain the following result:

\begin{thm} 
{\sl (a) Periods are divided into disjoint set of equivalent classes under the action of the
partial differential operators $\mathcal D_i$ ($i=1,2,\cdots, 5$). 

(b) Equivalently, the Picard--Fuchs differential
equation decomposes to the product of the Picard--Fuchs differential equations
corresponding to equivalent classes. 

(c) There is a one-to-one correspondence between the set of monomial classes of
under the action of $\hat G$-invariant monomials and the equivalence classes
of periods under the action of ${\mathcal D}_i$ operators.}
\end{thm}

\section{The monomial--motive correspondence: Examples}\label{monomialFermatcorrespondence}

In this section, we will formulate our main findings, namely,
a one-to-one correspondence between motives and monomials
for a mirror pairs of Calabi--Yau hypersurfaces in weighted projective $4$-spaces. This
correspondence is established at the Fermat points in the moduli spaces. 

First we illustrate the correspondence by two examples.  A monomial $X_1^{v_1}X_2^{v_2}X_3^{v_3}X_4^{v_4}X_5^{v_5}$
is represented by the exponent ${\mathbf v}=(v_1,v_2,v_3,v_4,v_5)$ and 
$(0,0,0,0,0)$ represents the constant.

\begin{expl} {\rm Consider the quintic hypersurfaces in the projective $\BP^4$. 
$$X_1^5+X_2^5+X_3^5+X_4^5+X_5^5-5\psi X_1X_2X_3X_4X_5=0.$$
So it has the weight $Q=(1,1,1,1,1)$. There are six classes of monomials (cf. Candelas et al. [CORV00, CORV03]).
Here multiplicity is the number of permutations.

\begin{table}
\begin{center}
\begin{tabular}{|c|c|c|c|c|}\hline
\#  & Monomial ${\mathbf v}$ & deg $w({\mathbf v})$& mult & length of orbit\\ \hline
1 & $(0,0,0,0,0)$& $0$ & $1$ & $5$ \\ \hline
2 & $(4,1,0,0,0)$& $5$ & $20$ & $5$ \\ \hline 
3 & $(3,2,0,0,0)$& $5$ & $20$ & $5$ \\ \hline
4 & $(3,1,1,0,0)$& $5$ & $30$ & $5$ \\ \hline
5 & $(2,2,1,0,0)$& $5$ & $30$ & $5$ \\ \hline
6 & $(4,0,3,2,1)$& $10$  &$24$ & $1$ \\ \hline
\end{tabular}
\vskip 0.2cm
\caption{Monomials for the quintic Calabi--Yau threefold} 
\end{center}
\end{table}

In the set of monomials, the $\hat G$-invariant monomials are
$$X_i^5\, (1\leq i\leq 5)\quad\mbox{and}\quad X_1X_2X_3X_4X_5.$$
The (repeated) multiplication by the $\hat G$-invariant monomial $X_1X_2\cdots X_5=(1,1,1,1,1)$ to each monomial
and reducing each component modulo $5$ establishes the monomial--motive correspondence.
The correspondences are:
$$[\ba]=[a_1,a_2,a_3,a_4,a_5]\mapsto \mbox{class of }\quad 
X_1^{a_1/5}X_2^{a_2/5}X_3^{a_3/5}X_4^{a_4/5}X_5^{a_5/5}$$
and
$${\mathbf v}=(v_1,v_2,v_3,v_4,v_5)\mapsto [{\mathbf v}+t(1,1,1,1,1)]\in
\frak{A}(Q)\quad\mbox{for $t\in(\Z/5\Z)^{\times}$}.$$
For instance, we have
 
$$(0,0,0,0,0)\mapsto (1,1,1,1,1)\leftrightarrow [1,1,1,1,1]$$
$$(4,1,0,0,0)\mapsto (1,3,2,2,2)\mapsto (3,4,1,1,1)\leftrightarrow [1,1,1,3,4]$$
$$(3,2,0,0,0)\mapsto (4,3,1,1,1) \leftrightarrow [1,1,1,3,4]$$
$$(3,1,1,0,0)\mapsto (4,2,2,1,1)\leftrightarrow [1,1,2,2,4]$$
$$(2,2,1,0,0)\mapsto (3,3,2,1,1)\mapsto (1,1,4,2,2)\leftrightarrow [1,1,2,2,4]$$

The length of orbit of each monomials of these actions is listed in the last column of the table.
 
The monomial $(4,0,3,2,1)$ only contributes at a conifold point, and hence it is not
realizable at the Fermat point.

\begin{table}
\begin{center}
\begin{tabular}{|c|c|c|c||c|c|c|}\hline \hline
\# & monomial ${\mathbf v}$ & deg $R([{\mathbf v}],t)$ & $\lambda_{\mathbf v}$ & motive & dim & mult \\ \hline
1 & (0,0,0,0,0) & 4 & 1 & [1,1,1,1,1] & 4 & 1 \\ \hline
2 & (4,1,0,0,0) & 4 & 20 & [1,1,1,3,4] & 4 & 20 \\
3 & (3,2,0,0,0) &  &   & \quad       &  &  \\ \hline
4 & (3,1,1,0,0) & 4 & 30 & [1,1,2,2,4] & 4 & 30 \\
5 & (2,2,1,0,0) &   &    &             &   &   \\ \hline
6 & (4,0,3,2,1) & 0 & 24 & ?          & ? &   \\ \hline
\end{tabular}
\vskip 0.2cm
\caption{Monomial--motive correspondence for the quintic Calabi--Yau threefold}
\end{center}
\end{table} 

The new notation introduced in TABLE 7 and TABLE 8 ought to be explained:
For a monomial ${\mathbf v}$, $R([{\mathbf v}],t)$ denotes the polynomial factor 
corresponding to the orbit $[{\mathbf v}]$ and $\lambda_{\mathbf v}$ its degree. (For details, see below.) 
\medskip

In Section 4, we observed that the motive corresponding to the weight $[1,1,1,1,1]$ is 
the only motive that is invariant under the mirror symmetry operation by the group 
$\hat G$. Also in Section 6, we showed in
the toric geometric setting, that the only the origin is left invariant under the mirror map.  
The monomial--motive correspondence is consistent with these two facts. } 
\end{expl}

\begin{expl} {\rm Consider the octic Calabi--Yau hypersurfaces in the weighted projective $4$-space
$\BP^4(Q)$ with $Q=(1,1,2,2,2)$ with two deformation parameters: 
$$X_1^8+X_2^8+X_3^4+X_4^4+X_5^5-2\phi X_1^4X_2^4-8\psi X_1X_2X_3X_4X_5=0$$
in $\BP^4(Q)$.  In this case, there are $15$ classes of monomials (see
Kadir [Ka04, Ka05]).

\begin{table}
\begin{center}
\begin{tabular}{|c|c|c|c||c|c|c|c|}\hline \hline
\# & monomial ${\mathbf v}$ & deg $w({\mathbf v})$  & mult & \# & monomial & deg & mult \\ \hline
1 & (0,0,0,0,0) & 0 & 1  & 9 & (0,4,0,3,3) & 16 & 3 \\ \hline
2 & (0,2,1,1,1) & 8 & 2  & 10 & (4,0,1,1,0) & 8 & 3 \\ \hline
3 & (6,2,0,0,0) & 8 & 1 & 11 & (2,0,3,0,0) & 8 & 6 \\ \hline
4 & (0,0,0,2,2) & 8 & 3 & 12 & (6,0,1,0,0) & 8 & 6 \\ \hline
5 & (2,0,1,3,3)& 16 & 6  & 13 & (0,0,3,1,0) & 8 & 6 \\ \hline
6 & (4,0,2,0,0) & 8 & 3 & 14 & (2,0,2,1,0) & 8 & 6 \\ \hline
7 & (0,0,2,1,1) & 8 & 3 & 15 & (4,0,3,2,1) & 16 & 6 \\ \hline
8 & (2,2,1,1,0) & 8 & 3 &    &             &   & \\ \hline
\end{tabular}
\vskip 0.2cm
\caption{Monomials for the octic Calabi--Yau threefold}
\end{center}
\end{table}

Here again the monomial $(4,0,3,2,1)$ only contributes at a conifold point, and hence
is not realizable at the Fermat point.

In this case, the $\hat G$-invariant monomials are
$$X_1^8, X_2^8, X_3^4, X_4^4, X_5^4$$
and
$$X_1^4X_2^4\quad\mbox{and}\quad X_1X_2X_3X_4X_5.$$

Given a monomial ${\mathbf v}=(v_1,v_2,v_3,v_4,v_5)$, adding 
$(1,1,1,1,1)$ or $(4,4,0,0,0)$ repeatedly
and then reducing the first two components modulo $8$ and the last three components
by modulo $4$ yield the monomial--motive correspondence.
The correspondences are given as follows with $Q=(q_1,q_2,q_3,q_4,q_5)=(1,1,2,2,2)$.
$$[\ba]=[a_1,a_2,a_3,a_4,a_5]\mapsto \mbox{class of}\,
\,X_1^{a_1/q_1}X_2^{a_2/q_2}X_3^{a_3/q_3}X_4^{a_4/q_4}X_5^{a_5/q_5}$$
and
$${\mathbf v}=(v_1,v_2,v_3,v_4,v_5)\mapsto [{\mathbf v}+\prod (1,1,2,2,2)^t (4,4,0,0,0)^s]
\in\frak{A}(Q)$$
for some $t,s\in \N$.
Here are some examples. First we describe motives to monomials correspondence.

$$[1,1,2,2,2]\mapsto X_1X_2X_3^{2/2}X_4^{2/2}X_5^{2/2}=\prod_{i=1}^5 X_i\mapsto (0,0,0,0,0)$$
$$[2,2,4,4,4]\mapsto X_1^2X_2^2X_3^{4/2}X_4^{4/2}X_5^{4/2}=(\prod_{i=1}^5 X_i)^2
\mapsto (0,0,0,0,0)$$
$$[5,1,2,4,4]\mapsto X_1^5X_2X_3^{2/2}X_4^{4/2}X_5^{4/2}=\prod_{i=1}^5 X_i(X_1^4X_4X_5)
\mapsto (4,0,0,1,1)$$
$$[3,1,2,4,6]\mapsto X_1^3X_2X_3X_4^2X_5^3=(\prod_{i=1}^5 X_i)(X_1^2X_4X_5^2)
\mapsto (2,0,0,1,2)$$
\medskip
Conversely, here are examples of monomials to motives correspondence.

$$(6,2,0,0,0)\mapsto [7,3,2,2,2]$$
$$(3,5,2,1,1)\mapsto (7,9,2,1,1)\mapsto (7,1,2,1,1)\mapsto [7,1,4,2,2]$$
\vskip 0.5cm

\begin{table}
\begin{center}
\begin{tabular}{|c|c|c|c||c|c|c|}\hline
\#& Monomial ${\mathbf v}$ &deg $R([{\mathbf v}],t)$ &$\lambda_{\mathbf v}$  &motive &dim&mult\\\hline
1 & $(0,0,0,0,0)$&$6$&$1$&$[1,1,2,2,2]$&$4$&$1$\\ &&&&$[2,2,4,4,4]$&$2$&$1$\\\hline
2 & $(0,2,1,1,1)$&$4$&$2$&$[3,1,4,4,4]$&$4$&$1$\\ &&&&$[6,4,2,2,2]$&$2$&$2$\\\hline
3 & $(6,2,0,0,0)$&$4$&$1$&$[7,3,2,2,2]$&$4$&$1$\\\hline
4 & $(0,0,0,2,2)$&$4$&$3$&$[3,3,2,2,6]$&$4$&$3$\\\hline
5 & $(2,0,1,3,3)$&$2$&$6$&$[4,2,2,2,6]$&$2$&$6$\\\hline
6 & $(4,0,2,0,0)$&$4$&$3$&$[5,1,2,2,6]$&$4$&$3$\\\hline
7 & $(0,0,2,1,1)$&$3$&$3$&$[4,4,2,2,4]$&$2$&$3$\\
8 & $(2,2,1,1,0)$&$3$&$3$&$[3,3,2,2,4]$&$4$&$3$\\\hline
9 & $(0,4,0,3,3)$&$4$&$3$&$[6,2,2,2,4]$&$2$&$6$\\
10 & $(4,0,1,1,0)$&$4$&$3$&$[5,1,2,2,4]$&$4$&$3$\\\hline
11 & $(2,0,3,0,0)$&$3$&$6$&$[4,2,2,4,4]$&$2$&$6$\\
12 & $(6,0,1,0,0)$&$3$&$6$&$[7,1,2,2,4]$&$4$&$6$\\\hline
13 & $(0,0,3,1,0)$&$2$&$6$&$[2,2,2,4,6]$&$2$&$6$\\\hline
14 & $(2,0,2,1,0)$&$4$&$6$&$[3,1,2,4,6]$&$2$&$12$\\\hline
15 & $(4,0,3,2,1)$&$0$&$6$&?&?&?\\\hline
\end{tabular}
\vskip 0.2cm
\caption{Monomial--motive correspondence for the octic Calabi--Yau threefold}
\end{center}
\end{table}
}
\end{expl}

In the above examples, the monomial--motive correspondences are obtained by
combinatorial matching of monomials and motives.
In the rest of this section we will establish this
correspondence  mathematically.  
\smallskip

First we re-access the situation.  Let $\mathcal R=\C[X_1,X_2,X_3,X_4,X_5]$ be the weighted polynomial ring of
the family of Calabi--Yau hypersurfaces in weighted projective $4$-space. 
 Let $\mathcal S$ denote the set of all
monomials ${\mathbf X}^{\mathbf v}=X_1^{v_1}X_2^{v_2}\cdots X_5^{v_5}$.  
Then $\mathcal S$ has the group of automorphisms consisting of
$\hat G$ and permutations of weighted coordinates.  Under the action of automorphisms,
$\mathcal S$ can be decomposed into the disjoint sets of orbits of monomials, which we
call a monomial class. A monomial class consists of monomials ${\mathbf X}^{\mathbf v}$
and its multiples by $\hat G$-invariant monomials. We may employ the notation 
$([v_1],[v_2],[v_3],[v_4],[v_5])$for the monomial class.

Since periods are in one-to-one correspondence with monomials, the decomposition
of monomials to the set of orbits under the action of automorphisms 
is carried over to periods.  On the set of periods, there are the differential
operators $\mathcal{D}_i$ described in Section 7. Under these differential
operators, periods decompose into the product of subperiods, which are
in one-to-one correspondence with the disjoint sets of orbits of monomials.    

\begin{thm} {\sl For a Calabi--Yau threefold $X$ associated to an admissible
pair $<m,Q>$ with $Q=(q_1,q_2,q_3,q_4,q_5)$, there is a one-to-one correspondence 
at the Fermat point between 
monomial classes and Fermat motives. The correspondences are given as follows:
A motive $A=[\ba]=[a_1,a_2,a_3,a_4,a_5]$ with $\ba\in\frak{A}(Q)$
corresponds to the monomial class 
$X_1^{a_1/q_1}X_2^{a_2/q_2}X_3^{a_3/q_3}X_4^{a_4/q_4}X_5^{a_5/q_5}$.

Conversely, 
a monomial class $X_1^{v_1}X_2^{v_2}X_3^{v_3}X_4^{v_4}X_5^{v_5}$ with
$([v_1],[v_2],[v_3],[v_4],[v_5])$ representing the class of monomials, corresponds to the motive 
$A=[[v_1]q_1,[v_2]q_2,[v_3]q_3,[v_4]q_4,[v_5]q_5]\in \frak{A}(Q)$.

(Note that the constant monomial class contain the monomial $(1,1,1,1,1)$ which corresponds
to the weight motive $\cM_Q=[q_1,q_2,q_3,q_4,q_5]$.)}
\end{thm}

\section{Proof of the monomial--motive correspondence}\label{proofcorrespondence}

We will prove Theorem 8.1 by comparing the two expressions for the congruence zeta-function of 
$X$ over $\F_q$ obtained in two different ways.  One is based on Weil's method and involves
Fermat motives. The other is based on Dwork's method and involves monomials. We will
show that the two approaches reconcile at the Fermat point giving the same congruence
zeta-function.  

Let $X$ be a Calabi--Yau threefold defined over $\Q$. 
To compute the number of $\F_q$-rational points on
the reduction $X_{\F_q}$ of $X$ to the finite field $\F_q$ of characteristic $p$.
(We need to choose a ``good prime'' $p$, for our examples, good primes
are those primes $p\nmid m$.) Denote by $\# X(\F_q)$ the number of $\F_q$-rational points on $X_{\F_q}$.
The congruence zeta function of $X_{\F_q}$ is defined by concocting the 
numbers $\# X_{\F_{q^r}}$ for all integers $r\geq 1$: 
$$Z(X_{\F_q},t)=\exp\left(\sum_{r=1}^{\infty}\# X_{\F_{q^r}}\frac{t^r}{r}\right).$$
Let $\bar{\F}_q$ denote the algebraic closure of $\F_q$ and write
$\bar X=X\otimes_{\F_q}\bar{\F}_q$. Choose a prime $\ell\neq p$. Then
by the work developed by Grothendieck and reached the summit by Deligne [D74],
$Z(X_{\F_q},t)$ has the interpretation in terms of $\ell$-adic etale cohomology groups 
$H^i_{et}(\bar X, \Q_{\ell})$ as
$$Z(X_{\F_q},t)=\prod_{i=0}^6 P_i(X_{\F_q}, t)^{(-1)^{i-1}}$$
where $$P_i(X_{\F_q}, t):=det(1-t\,Frob_q\,|\, H^i(\bar X,\,Q_{\ell}))\in 1+\Z[t]$$
is an integral polynomial of degree equal to the $i$-th Betti number $B_i(X)$, and
it satisfies the Riemann Hypothesis, namely, its reciprocal roots are algebraic
integers with absolute value $q^{i/2}$.
\medskip

For diagonal hypersurfaces, the congruence zeta-functions were explicitly described
in Weil [We49]. The number $\# X(\F_q)$ were computed using Jacobi sums and Gauss sums.
For our Calabi--Yau orbifolds associated to admissible pairs $<m,Q>$, the congruence
zeta-functions are of the form
$$Z(X_{\F_q},t)=\frac{P_3(X_{\F_q},t)}{(1-t)(1-qt)^{h^{1,1}(X)}(1-q^2t)^{h^{2,2}(X)}(1-q^3t)}$$
only the middle cohomology group $H^3$ (and hence $P_3$) contains essential information.

\begin{prop} {(\rm Yui [Y05])} {\sl Let $X$ be a Calabi--Yau orbifold
corresponding to an admissible pair $<m,Q>$ with $Q=(q_1,q_2,q_3,q_4,q_5)$. 
Let $q$ be a power of prime $p$ such that $q\equiv 1\pmod m$.
Then
$$P_3(X_{\F_q},t)=\prod_{\ba\in\frakA(Q)} (1-J(\ba)\,t)$$
with 
$$\frakA(Q)=\{\ba\in(\Z/m\Z)^5\,|\, a_i\neq 0,\, \sum_{i=1}^5 a_i\equiv 0\pmod m\,\}.$$
Here $J(\ba)$ is the Jacobi sum defined as
$$J(\ba)=(-1)^3\sum \chi(x_1)^{a_1}\chi(x_2)^{a_2}\chi(x_3)^{a_3}\chi(x_4)^{a_4}\chi(x_5)^{a_5}$$
where $\chi:\F_q\to \mu_m$ is a multiplicative character of order $m$,
and the sum runs over $(x_1,x_2,\cdots,x_5)\in (\F_q^*)^5$ such that $\sum_{i=1}^5x_i=0\in\F_q$.  

Furthermore, the Jacobi sum is expressed as the product of Gauss sums as follows:
$$J(\ba)=(-1)^3\frac{1}{q}G(\chi^{a_1})G(\chi^{a_2})\cdots G(\chi^{a_5})$$
where $G(\chi)$ is the Gauss sum defined by
$$G(\chi)=G(\chi,\psi):=\sum_{x\in\F_q^*} \chi(x)\psi(x)$$
with $\psi$ an additive character of $\F_q$.

(If we use a different normalization by putting minus sign in front of the sum
in the definition of Gauss sums, then $G(\chi_1,\phi)=1$ for the trivial character $\chi_1$.
In this case, the minus sign in front of $J(\ba)$ can be dropped.  There are two types
of normalizations in literature.)

Furthermore, $P_3(X,t)$ factors into the product:
$$P_3(X_{\F_q},t)=\prod_{A\in O(\frakA(Q))} P_3(\cM_A,t)$$
where
$$P_3(\cM_A,t)=\prod_{\ba\in A}(1-J(\ba)\,t)\in 1+\Z[t]\quad\mbox{with deg\, $P_3(\cM_A,t)=B_3(\cM_A)$}.$$}
\end{prop}

By Weil [We49], we know that Gauss sums are algebraic integers with absolute
value $q^{1/2}$ in the field $\Q(\zeta_m,\zeta_q)$ over $\Q$ where $\zeta_m$ (resp.
$\zeta_q$) is a primitive $m$-th (resp. $q$-th) root of unity.
Therefore, Jacobi sums $J(\ba)$ are algebraic integers in $\bL=\Q(\zeta_m)$
with absolute value $q^{3/2}$. 
Weil's formula for $|X(\F_q)|$ in terms of Gauss and Jacobi sums.

\begin{lem}{\rm (Weil's formula)} {\sl Let $X$ be a Calabi--Yau threefold corresponding to an admissible
pair $<m,Q>$. Then the number of $\F_q$-rational point on $X$ is given by the
formula:
$$\#X(\F_q)=1+q+q^2+q^3-\sum_{\ba\in\frakA(Q)} J(\ba)=1+q+q^2+q^3+\frac{1}{q}\sum_{\ba\in\frakA(Q)}\prod_{i=1}^5
G(\chi^{a_i}).$$}
\end{lem}

We may further factor the last term into sums involving Fermat motives as
$$N_{motive}:=\sum_{\ba\in\frakA(Q)} J(\ba)=\sum_{A\in O(\frakA(Q))}\sum_{\ba\in A} J(\ba).$$

Now we consider a mirror partner $\hat X$ of a Calabi--Yau orbifold $X$ corresponding to
$<m,Q>$, and its congruence zeta-function at the Fermat point.

\begin{thm} {\sl Let $(X, \hat X)$ be a mirror pair of Calabi--Yau orbifolds over $\Q$ corresponding to
an admissible pair $<m,Q>$ constructed in Section 4.  Let $X_{\F_q}$ (resp. $\hat X_{\F_q}$) be the
reduction of $X$ (resp. $\hat X$) to the finite fields $\F_q$ of characteristic $p$ (where
$p$ is a good prime). Then the congruence zeta-function of $\hat X$ at the
Fermat point is given as follows:
$$Z(\hat X_{\F_q},t)=\frac{P_3(\hat X_{\F_q}, t)}{(1-t)(1-qt)^{h^{2,1}}(1-q^2t)^{h^{1,2}}(1-q^3t)}$$
with
$$P_3(\hat X_{\F_q},t)=\prod_{A} P_3(\cM_A,t)$$ 
where the product runs over all $A$ such that corresponding Fermat motives $\cM_A$ are
invariant under the mirror map (i.e., the action of $\hat G$).}
\end{thm}

\medskip

There is another method of calculating the number of $\F_q$-rational points on
the family of Calabi--Yau hypersurfaces.  This method was developed by Dwork [D60]. 
The number of $\F_q$-rational point is computed using $p$-adic analytic method, namely, $p$-adic
Gamma function, the Dwork character, and the Gross--Koblitz formula. This method works for
more general hypersurfaces, e.g., deformations of diagonal hypersurfaces.  

For this, we review the Dwork character, Gauss sums, the Gross--Koblitz formula from Candelas et al. [CORV00]
(tailored especially for physicists). 

Let $\F_p$ be the prime field of $p$ elements and let $\Z_p^*\simeq \mu_{p-1}$
be the group of $(p-1)$-th roots of unity. The Teichm\"uller character
$\omega_p\,:\, \F_p^*\to \Z_p^*$ is a multiplicative character of $\F_p^*$ of order $p-1$,
so that $\omega_p(x)\equiv x\pmod p$.  Let $\pi$ be an element in the algebraic closure $\bar{\Q}_p$ such that
$\pi^{p-1}=-p$, and let $F$ be the function defined by
$F(x)=exp(X+X^p/p).$ 
The Dwork character $\Theta$ is an additive character of $\F_p$ defined by
$$\Theta(x)=\Theta_{p^0}(x)=F(\pi \omega_p(x)).$$
Then we can define a Gauss sum 
$$G_n=G(\Theta, \omega_p^n)=\sum_{x\in\F_p^*} \Theta(x) \omega_p^n (x).$$
We have the inversion formula
$$\Theta(x)=\frac{1}{p-1}\sum_{k=0}^{p-2}G_{-k} \omega_p^k(x).$$

For a finite extension $\F_q$ ($q=p^r)$)  of $\F_p$, 
the Teichm\"uller character $\omega_q$ is defined as a multiplicative character of
$\F_q^* \to \mu_{q_1}$ of the maximal order $q-1$, and
the Dwork character $\Theta_r$ is an additive
character of $\F_q$ defined by composing $\Theta_{p^0}$ with the trace
$\mbox{tr}:\F_q\to \F_p,\,\,\mbox{tr}(x)=x+x^p+x^{p^2}+\cdots +x^{p^{r-1}}.$
That is,
$$\Theta_r(x)=\prod_{i=0}^{r-1}\Theta_{p^i}(x).$$
A Gauss sum is defined by
$$G_{r,n}=G(\Theta_r,  \omega_q^n)=\sum_{x\in\F_q^*} \Theta_r(x) \omega_q^n(x).$$

\begin{lem} (The Gross--Koblitz formula) {\sl Let $\omega_p$ be the Teichm\"uller
character of $\F_p^*$. Let $\Gamma_p$ be the $p$-adic Gamma function defined by
$$\Gamma_p(n):=(-1)^n\prod_{1\leq i<n, p\nmid i} i$$
taking values in the $p$-adic integer ring $\Z_p$. 
Then the Gauss sum $G(\Theta,\omega_p)$ is expressed by the Gross--Koblitz
formula for the case $p-1\nmid n$:
$$G_n=p(p)^{-\langle\frac{n}{p-i}\rangle} \Gamma_p\left(1-\left\langle\frac{n}{p-1}\right\rangle\right).$$
Here $\langle x \rangle$ denotes the fractional part for $x\in\R$.

For a finite extension $\F_q$ ($q=p^r$), the Gross--Koblitz formula for
the case $q-1\nmid n$ is given by
$$G_{r,n}=G(\Theta_r,\omega_q^n)=(-1)^{r+1}q\pi^{-S(n)}\prod_{i=0}^{r-1} \Gamma_p\left(1-\left\langle\frac{p^in}{q-1}\right\rangle\right)$$
in $\Q_p(\pi)$ where $S(n)$ denotes the sum of the $p$-adic digits of $n$.}
\end{lem}

\begin{lem} (Dwork's formula) {\sl Let $X$ be a variety defined as the zero set of
a polynomial $P({\mathbf X})\in \F_q[X_1,X_2,\cdots, X_5]$. Then
$$\sum_{y\in \F_q}\Theta(yP(X))=\begin{cases} 0 & \mbox{if $P(X)\neq 0$}\\
                                           q &\mbox{if $P(X)=0$}\end{cases}$$
Hence we have
$$q\# X(\F_q)=\sum_{X_i\in \F_q}\sum_{y\in \F_q}\Theta(yP({\mathbf X})).$$}
\end{lem}

The congruence zeta-function $Z(X_{\F_q},t)$ of $X_{\F_q}$ is defined same as above
by counting the number $\# X_{\F_q}$ for each $r\geq 1$. The Dwork theory has
established the rationality of $Z(X_{\F_q},t)$ ([Dw60].) 
We can compute the congruence zeta-function for our mirror pairs of Calabi--Yau
threefolds.

Let $<m,Q>$ be an admissible pair, and let $\Delta(Q)$ be a reflexive polytope. 
Let $f({\mathbf c},{\mathbf X})=0\in\BP^4(Q)$ be a defining equation of degree $m$ over $\Q$ with 
deformation parameter ${\mathbf c}$.  Let $X_{\Delta(Q)}$ denote a family of
Calabi--Yau threefolds.
Let ${\mathbf v}=(v_1,v_2,v_3,v_4,v_5)$ denote
a monomial $X_1^{v_1}X_2^{v_2}\cdots X_5^{v_5}$ in the weighted polynomial coordinate ring
$\C[X_1,X_2,\cdots, X_5]$ of $X_{\Delta(Q)}$, and let $\lambda_{\mathbf v}$ stands for
its multiplicity, i.e., the number of weighted permutations. 
Let ${\mathcal S}$ be the set of all monomials. Then we know that ${\mathcal S}$ decomposes
into the disjoint sets $[{\mathbf v}]$ of orbits of monomials ${\mathbf v}$ under the action of the
automorphism group. Denote by $O({\mathcal S})$ the set of orbits.

\begin{thm} {\sl Let $<m,Q>$  be an admissible pair, and let $X_{\Delta}$ be a family of Calabi--Yau threefolds. 
For a good prime $p$, let $X_{\Delta/\F_q}$ denote the reduction of $X_{\Delta}$ over the finite field
$\F_q$. 
Then the congruence zeta-function of $X_{\Delta/\F_q}$ is a rational function
over $\Z$ and is given by
$$Z(X_{\Delta/\F_q},t)=\frac{R(X_{\Delta/\F_q},t)}{(1-t)(1-qt)^{h^{1,1}(X_{\Delta})}
(1-q^2t)^{h^{2,2}(X_{\Delta})}(1-q^3t)}.$$
Here $R(X_{\Delta/\F_q},t)$ is an integral polynomial of degree $B_3(X_{\Delta})$ and
has the form:
$$R(X_{\Delta(Q)/\F_q},t)=\prod_{[{\mathbf v}]\in O({\mathcal S})} R([{\mathbf v}],t)
=R({\mathbf 0},t)\prod_{[{\mathbf v}]\neq [{\mathbf 0}]} R({[\mathbf v}],t).$$
Furthermore, $R([{\mathbf v}],t)$ is an integral polynomial (not necessarily
irreducible), and
$$\sum_{[{\mathbf v}]\in O({\mathcal S})} \lambda_{\mathbf v} deg(R)=B_3(X_{\Delta(Q)})$$ 
where $\lambda_{\mathbf v}$ denotes the multiplicity of ${\mathbf v}$.  }\label{Fermatpointtheorem}
\end{thm} 

Now consider a mirror partner $X_{\Delta^*(Q)}$ and its congruence zeta-function.

\begin{thm} {\sl Let $(X_{\Delta(Q)}, X_{\Delta^*(Q)})$ be a mirror pair 
of Calabi--Yau threefolds defined over $\Q$ corresponding to a pair of reflexive polytopes
$(\Delta(Q), \Delta^*(Q))$ constructed in Section 6. Let $X_{\Delta(Q)/\F_q}$
(resp. $X_{\Delta^*(Q)/\F_q}$) be the reduction of $X_{\Delta(Q)}$ (resp.
$X_{\Delta^*(Q)}$) to the finite field $\F_q$ of characteristic $p$ (where
$p$ is a good prime). Then the congruence zeta-function of $X_{\Delta^*(Q)/\F_q}$
at the Fermat point is given as follows:
$$Z(X_{\Delta^*(Q)/\F_q},t)=\frac{R([{\mathbf 0}],t)}
{(1-t)(1-qt)^{h^{2,1}(X_{\Delta(Q)})}(1-q^2t)^{h^{1,2}(X_{\Delta(Q)})}(1-q^3t)}$$ 
where $R([{\mathbf 0}],t)$ is defined as in Theorem 9.6.}
\end{thm}

We compare the expression for the congruence zeta-function for $X$ over $\F_q$
obtained from Weil's method and that from Dwork method at the Fermat point in the 
moduli space (i.e.,  at the point where all parameters are set to zero). 

\begin{thm}{\sl Let $<m,Q>$ be an admissible pair. Let $(X, \hat X)$ be a mirror
pair of Calabi--Yau threefolds. The for a good prime $p$, the two expressions
for the denominator of the congruence zeta-function for $X$ coincide at the
Fermat point. That is,
$$P_3(X_{\F_q},t)=\prod_{A\in O(\frakA(Q))} P_3(\cM_A,t)^{dim (A)}=R(X_{\Delta(Q)/\F_q},t)
=\prod_{[{\mathbf v}]\in O({\mathcal S})} R({\mathbf v},t)^{\lambda({\mathbf v)}}.$$

In particular, the factor $R([{\mathbf 0}],t)$ coincides with the factor
$\prod_{A} P_3(\cM_A,t)$ where $A$ runs over the $\hat G$-invariant motives.
Moreover, this is the only factor present in the denominator of the congruence
zeta function of $X$ and $\hat X$.
}
\end{thm}

For the calculation of the polynomials $P_3(\cM_A,t)$, we have the following
results for the cases when the polynomials become powers of linear polynomials over $\Q$,
and also when the polynomials are irreducible over $\Q$.

\begin{prop} {\rm (Yui[05])} {\sl Fix an admissible pair $<m,Q>$, and
let $\cM_A$ be a Fermat motive of $\ba\in\frakA(Q)$.  
Let $p$ be a good prime and let $f$ be the order of $p$ mod $m$, i.e.,
$f$ is the smallest positive integer such that $p^f\equiv 1\pmod m$. 

(a) If $f$ is even and $p^{f/2}+1\equiv 0\pmod m$, then putting $q=p^f$, we have
$$P_3(\cM_A,t)=(1-q^{3/2}t)^{\varphi(m)}.$$

(b) Let $\cM_Q$ be the motive corresponding to the weight $Q=(q_1,q_2,q_3,q_4,q_5)$.
Then for $p\equiv 1\pmod m$, $P_{\cM_Q}(t)$ is always irreducible over $\Q$.}
\end{prop}

\begin{expl} {\rm For $m=5$ and $Q=(1,1,1,1,1)$, the polynomials
$P_3(\cM_A,t)$ are given as follows.  By the (Newton) slopes, we mean the (normalized) $p$-adic
order of the reciprocal roots of $P_3(\cM_A,t)$, so that they take rational values in $[0, 3]$
with multiplicity $(*)$.
\smallskip

(a) When $p\equiv 1\pmod 5$ (e.g., $p=11$), we have
$$P_3([1,1,1,1,1],t)=1+89t+3^3\cdot 11\cdot 13t^2+11^3\cdot 89t^3+11^6t^4\quad\mbox{with
slopes $0,1,2,3$}$$
$$P_3([1,1,1,3,4],t)=1-11t-3^211^2t^2-11^4t^3+11^6t^4=P_3([1,1,2,2,4],t)\quad\mbox{with
slopes $1 (2), 2 (2)$}.$$
Both polynomials are irreducible over $\Q$.

(b) When $p\not\equiv 1\pmod 5$, then $p\equiv 2,3$ or $4 \pmod 5$. 
If $p\equiv 2, 3\pmod 5$, then $p^4\equiv 1\pmod 5$ and $p^2\equiv -1\pmod 5$. In this
case, for every motive $\cM_A$, we have 
$$P_3(\cM_A,t)=(1-p^{4\cdot\frac{3}{2}}t)^4\quad\mbox{with slopes $3/2 (4)$}.$$

If $p\equiv 4\pmod 5$, then $p^2\equiv 1\pmod 5$ and $p\equiv -1\pmod 5$. In this case,
for every motive $\cM_A$, we have
$$P_3(\cM_A,t)=(1-p^{2\cdot \frac{3}{2}}t)^4\quad\mbox{with slopes $3/2 (4)$}.$$
\smallskip

On the other hand, using Dwork's method, the polynomials $R([{\mathbf 0}],t)$, $R([4,1,0,0,0],t)$
and $R([3,2,0,0,0],t)$ are computed by Candelas et al. [CORV00, CORV03] as follows.
For $p=11$,
$$R([{\mathbf 0}],t)=1+89t+3^3\cdot 13\cdot 11t^2+89\cdot 11^3t^3+11^6t^4\quad\mbox{with slopes
$0,1,2,3$}$$
and
$$R([4,1,0,0,0],t)=R([3,2,0,0,0],t)=1-11t-3^2 11^2t^2-11^4t^3+11^6t^4$$
with slopes $1 (2),2(2)$.}
\end{expl}

\begin{expl} {\sl For $m=8$ and $Q=(1,1,2,2,2)$, the polynomials $P_3(\cM_A,t)$ are given
as follows.

(a) When $p\equiv 1\pmod 8$ (e.g., $p=17$), we have
$$P_3([1,1,2,2,2],t)=1-2^23^2\cdot 5t+2\cdot 17\cdot 467t^2-2^23^2\cdot 5\cdot 17^3t^3+17^6t^4\quad
\mbox{with slopes $0,1,2,3$}$$
$$P_3([1,1,2,6,6],t)=(1+2\cdot 3\cdot 17t+17^3t^2)^2=P_3([1,1,4,4,6],t)$$
$$P_3([1,3,2,4,6],t)=(1-2\cdot 3\cdot 17t+17^3t^2)^2=P_3([1,3,4,4,4],t)$$
$$P_3([1,5,2,2,6],t)=(1-2\cdot 17t+17^3t^2)^2=P_3([1,5,2,4,4],t)=P_3([1,5,6,6,6],t)$$
all with slopes $1(2), 2(2)$. For all the remaining motives (see TABLE 4), 
$$P_3(\cM_A,t)=(1+2\cdot 17t+17^3t^2)^2\quad\mbox{with slopes $1(2), 2(2)$}.$$
In particular,
$$P_3([2,2,4,4,4],t)=(1+2\cdot 17t+17^3t^2)^2.$$

(b) When $p\not\equiv 1\pmod 8$, $p\equiv 3,5$ or $7\pmod 8$, and in all cases
$p^2\equiv 1\pmod 8$.

If  $p\equiv 7\pmod 8$, we have, for very motive $\cM_A$, 
$$P_3(\cM_A,t)=(1-p^3t)^4\quad\mbox{with slopes $3/2 (4)$}.$$

If $p\equiv 3\pmod 8$ (e.g., for $p=11$), we have 
$$P_3([1,1,2,2,2],t)=(1-2\cdot 7\cdot 11^2t+11^6t^2)^2=P_3([1,1,2,6,6],t)=P_3([1,1,4,4,6],t)$$
$$P_3([1,3,2,4,6],t)=(1+2\cdot 7\cdot 11^2t+11^6t^2)^2=P_3([1,3,4,4,4],t)$$
all with slopes $1(2), 2(2)$. For the remaining motives $\cM_A$, 
$$P_3(\cM_A,t)=(1\pm 11^3t)^4\quad\mbox{with slopes $3/2 (4)$}.$$
In particular,
$$P_3([2,2,4,4,4],t)=(1-11^3t^2)^4.$$

If $p\equiv 5\pmod 8$ (e.g., for $p=13$), we have
$$P_3([1,1,2,2,2],t)=(1+2\cdot 7\cdot 13\cdot 17t+13^6t^2)^2\quad\mbox{with slopes $1(2), 2(2)$}$$
$$P_3([1,1,2,6,6],t)=(1-13^3t)^4=P_3([1,1,4,4,6],t)=P_3([1,3,2,4,6],t)=P_3([1,3,4,4,4],t)$$
with slopes $3/2 (4)$. For the remaining motives $\cM_A$,
$$P_3(\cM_A,t)=(1\pm 2\cdot 5\cdot 13^2+13^6t^2)^2\quad\mbox{with slopes $1(2), 2(2)$}.$$
In particular,
$$P_3([2,2,4,4,4],t)=(1+2\cdot 5\cdot 13^2t+13^6t^2)^2.$$
\smallskip

Kadir [Ka04] computed the polynomials $R([{\mathbf v}],t)$ using Dwork's method, the results are shown in Table 10. The following notation is used:

\begin{center}
\begin{tabular}{|c|c|}\hline
Notation&Polynomial, for prime $p$\\\hline
$(a)_1$&$(1+at)$\\
$(a)_2$&$(1+at+p^3t^2)$\\
$(a,b)_4$&$(1+at+bt^2+ap^3t^3+p^6t^4)$\\
\hline
\end{tabular}
\end{center}

\begin{table}
{\small
\begin{tabular}{|c|c|c|c|}\hline
Monomial $\textbf{v}$ & $p=11$ & $p=13$ &$p=17$\\
\hline\hline
$(0,0,0,0,0)$ & $(0)_2(0,2.7.11^2)_4$ &  $(6.13)_2(0,-238.13)_4$ & $(-2.17)_2(180,934.17)_4$\\
$(0,2,1,1,1)$& $(0)_2[(-6.11)_2(6.11)_2]^{1/2}$ &$(0)_2(-6.13)_2$ & $(-2.17)_2(6.17)_2$\\
$(6,2,0,0,0)$ & $(0)_2^2$ & $(-4.13)_2(4.13)_2$& $(2.17)_2^2$ 
\\ \hline
$(0,0,0,2,2)$ & $(0,2.7.11^2)_4$& $(0)_2$ & $(-6.17)_2^2$ \\
$(2,0,1,3,3)$ & $(0)_2$ & $(-6.13)_2$&$(-2.17)_2$  \\
$(4,0,2,0,0)$ & $(0)_2^2$& $(-4.13)_2(4.13)_2$& $(2.17)_2^2$
\\ \hline
$(0,0,2,1,1)$&$(0)_2^{1/2}(0,2.7.11^2)_4^{1/2}$&$(0)_2(6.13)_2$&$(-6.17)_2^2(-2.17)_2$ \\
$(2,2,1,1,0)$&$(0)_2^{1/2}(0,2.7.11^2)_4^{1/2}$&                              &\\\hline
$(6,0,1,0,0)$&$(0)_2^{1/2}(11\sqrt{11})_1(-11\sqrt{11})_1$&$(-4.13)_2(4.13)_2(6.13)_2$&$(-2.17)_2^3$\\
$(2,0,3,0,0)$&$(0)_2^{1/2}(11\sqrt{11})_1(-11\sqrt{11})_1$&                                                    &\\\hline
$(0,4,0,3,3)$&$(0)_2^2$&$(-4.13)_2(4.13)_2\times$&$(-2.17)_2^2(2.17)_2^2$\\
$(4,0,1,1,0)$&$(0)_2^2$&$(-6.13)_2(6.13)_2$&\\\hline
$(0,0,3,1,0)$ & $(0)_2$  & $(-6.13)_2$ & $(-2.17)_2$\\
$(2,0,2,1,0)$ & $[(-6.11)_2(6.11)_2]^{1/2}$&$(0)_2$ & $(6.17)_2$\\
$(4,0,2,3,1)$ & $1$  & $1$&$1$ \\ \hline
\end{tabular}
}
\vskip 0.3cm
\caption{Polynomials for $p=11,13,17$}
\end{table}
}
\end{expl}
 
Alternatively, since the basic constituents of the congruence zeta-function are the
number of $\F_{q^r}$-rational points on $X$, we compare
the expression for the number of $\F_q$-rational points, $N_{motive}$ by Weil's method
and the number of $\F_q$-rational pints, $N_{mon}$ by Dwork's method.    


\medskip

\noindent {\bf Proof of Theorem 9.8.}
We now prove Theorem 9.8 by showing that $N_{motive}=N_{mon}.$ we  
now count the number of points on the variety excluding the contribution coming from any 
exceptional divisor required to smooth the ambient weighted projective space. 
The piece coming from the exceptional divisor needs to be added by hand using Weil's method, 
but emerges naturally from toric geometry using Dwork's method (namely toric geometry forbids 
the vanishing of coordinates associated to points lying on the same cone; see [CORV00, CORV03] and [Ka04, Ka05] 
for examples).  Using Dwork's method, as outlined above, we obtain the following expression for 
the case when we have an admissible pair $<m,Q>$ with $Q=(q_1,q_2,q_3,q_4,q_5)$ :

\begin{eqnarray*}
q\nu =\sum_{y\in\mathbb{F}_p}\sum_{\textbf{x}\in(\mathbb{F}_p)^5}\Theta\left(y\left(\sum_{i=1}^{5}x_i^{\frac{m}{q_i}}\right)\right) 
=q^5 + \sum_{y\in\mathbb{F}_p^\ast}\sum_{\textbf{x}\in(\mathbb{F}_p)^5}\Theta\left(y\left(\sum_{i=1}^{5}x_i^{\frac{m}{q_i}}\right)\right)+\nu_{exceptional}.\\\nonumber
\end{eqnarray*}
Dividing both sides by $q$ and bringing in Gauss sums, we get:
\begin{eqnarray*}
\nu
&=& q^4 + \frac{1}{q(q-1)^5}\sum_{y\in\mathbb{F}_p^\ast}\sum_{\textbf{x}\in(\mathbb{F}_p^\ast)^5}\prod_{i=1}^{5} G_{-s_i}\omega_q^{\frac{ms_i}{q_i}} (x_i)\omega_q^{\sum s_i} (y)
\\\nonumber
&&  -G_0\frac{1}{q(q-1)^4}\sum_{j=1}^5\sum_{y\in\mathbb{F}_p^\ast}\sum_{x_i\in\mathbb{F}_p^\ast\atop i\neq j}\prod_{i\neq j} G_{-s_i}\omega_q^{\frac{ms_i}{q_i}} (x_i)\omega_q^{\sum s_i} (y)
\\\nonumber
&&  +G_0^2\frac{1}{q(q-1)^3}\sum_{j, k}\sum_{y\in\mathbb{F}_p^\ast}\sum_{x_i\in\mathbb{F}_p^\ast\atop i\neq j}\prod_{i\not\in {\{j,k\}}} G_{-s_i}\omega_q^{\frac{ms_i}{q_i}} (x_i)\omega_q^{\sum s_i} (y)
\\\nonumber
&&  -G_0^3\frac{1}{q(q-1)^2}\sum_{j, k,l}\sum_{y\in\mathbb{F}_p^\ast}\sum_{x_i\in\mathbb{F}_p^\ast\atop i\neq j}\prod_{i\not\in {\{j,k,l\}}} G_{-s_i}\omega_q^{\frac{ms_i}{q_i}} (x_i)\omega_q^{\sum s_i} (y)
\\\nonumber
&&  +G_0^4\frac{1}{q(q-1)}\sum_{j, k,l,n}\sum_{y\in\mathbb{F}_p^\ast}\sum_{x_i\in\mathbb{F}_p^\ast\atop i\neq j}\prod_{i\not\in {\{j,k,l,n\}}} G_{-s_i}\omega_q^{\frac{ms_i}{q_i}} (x_i)\omega_q^{\sum s_i} (y)
\\\nonumber
&& -G_0^5\frac{q-1}{q}+\nu_{exceptional}.\
\end{eqnarray*}

Further, we pass onto the following expression for $\nu$:
\begin{eqnarray*}
\nu
&=& q^4 + \frac{q-1}{q}\sum_{s_i=1}^{q-2}\prod_{i=1}^{5} G_{-s_i}
 -G_0\frac{q-1}{q}\sum_{j=1}^5\sum_{s_i=1\atop i\neq j}^{q-2}\prod_{i\neq j} G_{-s_i}\\\nonumber
&& +G_0^2\frac{q-1}{q}\sum_{j, k}\sum_{s_i=1\atop i\not\in {\{j,k\}}}^{q-2}\prod_{i\not\in {\{j,k\}}} G_{-s_i}  -G_0^3\frac{q-1}{q}\sum_{j, k,l}\sum_{s_i=1\atop i\not\in{\{j,k,l\}}}^{q-2}\prod_{i\not\in {\{j,k,l\}}} G_{-s_i}\\\nonumber
&&  +G_0^4\frac{q-1}{q}\sum_{j, k,l,n}\sum_{s_i=1\atop i \not\in{\{j,k,l,n\}}}^{q-2}
\prod_{i\not\in {\{j,k,l,n\}}} G_{-s_i} -G_0^5\frac{q-1}{q}+\nu_{exceptional},\\\nonumber
\end{eqnarray*}
where passing from the former to the latter expression we had to impose the following conditions:
\begin{center}
\begin{math}
\begin{array}{ll}
q-1|\quad\frac{ms_i}{q_i}&\forall i=1,\ldots, 5\\
q-1|\quad\sum_{i=1}^{5}s_i& 
\end{array}
\end{math}
\end{center}

\noindent Hence introduce a vector of integers, $\textbf{v}^{\textbf{t}}=(v_1,v_2,v_3,v_4,v_5)$, such that:\\
\begin{center}
\begin{math}
\begin{array}{llll}
v_i(q-1)&=&\frac{s_i m}{q_i}&\forall i=1,\ldots, 5\\
v(q-1)&=&\sum_{i=1}^{5}s_i&
\end{array}
\end{math}
\end{center}
and of course $0\leq s_i\leq q-2, \quad i=1,\ldots,5$.
It is easy to see that we need to consider vectors $\textbf{v}$
such that:
\begin{eqnarray*}
\sum_{i=1}^{5}q_iv_i= mv ,\quad
0\leq q_iv_i\leq \bigg\lfloor\frac{m(q-2)}{q-1}\bigg\rfloor \quad(\forall i=1,\ldots, 5)\nonumber\\
\end{eqnarray*}
\begin{eqnarray*}
0\leq v\leq \bigg\lfloor\frac{5(q-2)}{q-1}\bigg\rfloor \quad(= 4\quad \text{for}\; q\geq 7)\nonumber\\
\end{eqnarray*}

Hence,

\begin{eqnarray*}
\nu
&=& q^4 + \frac{q-1}{q}\sum_{\textbf{v}\atop\prod v_i\neq 0}\prod_{i=1}^{5} G_{-\left(\frac{q-1}{m}v_iq_i\right)}+\nu_{exceptional}.
\\\nonumber
\end{eqnarray*}

The reason all monomials $\textbf{v}$ with zero entries are excluded is as follows: Consider a 
monomial $\textbf{v}$ with exactly $r$ zero entries. In the second equation for $\nu$, it will appear 
exactly ${r \choose 0}$ times in the first sum, and $(-1)^u{r \choose u}$ times in the $(u+1)$-th sum. 
Hence in total it appears $\sum_{u=0}^r (-1)^u{r \choose u} = 0$ times. 
Clearly this argument works in any dimension.

Finally we obtain:
\begin{eqnarray*}
N_{mon}
= \frac{1}{q}\sum_{\textbf{v}\atop\prod v_i\neq 0}\prod_{i=1}^{5} G_{-\left(\frac{q-1}{m}v_iq_i\right)}
= \frac{1}{q}\sum_{\textbf{v}\atop\prod v_i\neq 0}\prod_{i=1}^{5} G_{\left(\frac{q-1}{m}\right)(m-v_iq_i)}\;,
\end{eqnarray*}
where the last equality comes from the fact that $G_{r,n}=G_{r,q-1+n}$. If we define the following:
$$
q_iv_i= a_i^{\prime} ,\quad a_i= m-a_i^{\prime} ,\quad
\forall i=1,\ldots 5 $$
we see that
\begin{eqnarray*}
vm=\sum_{i=1}^5a_i=\sum_{i=1}^5(m-a_i^{\prime})=m(5-v^{\prime}) .\nonumber\\
\end{eqnarray*}
 It is clear that the vectors with components $a_i$ are precisely those in $\frakA(Q)$; hence
  \begin{eqnarray*}
N_{mon}&=&  \frac{1}{q}\sum_{\textbf{v}\atop\prod v_i\neq 0}\prod_{i=1}^{5} G_{\left(\frac{q-1}{m}a_i\right)}
\\\nonumber
\end{eqnarray*}
On the other hand, we know that the expression for the sum $\sum_{\ba\in\frakA(Q)} J(\ba)$ in terms of 
Weil's method is
$$\frac{1}{q}\sum_{\ba\in\frakA(Q)}\prod_{i=1}^5 G(\chi^{a_i}).$$ 
Therefore, we obtain: 
$$N_{motive}= N_{mon}.$$
This completes the proof of Theorem 9.8.
\medskip

As a consequence of the proof of Theorem 9.8, we can now establish the motive--monomial correspondence
at the Fermat point.

For a character $\ba=(a_1,a_2,a_3,a_4,a_5)\in\frakA(Q)$, there corresponds the monomial 
$$X_1^{a_1/q_1}X_2^{a_2/q_2}X_3^{a_3/q_3}X_4^{a_4/q_4}X_5^{a_5/q_5}.$$ 
Now the $(\Z/m/Z)^{\times}$-orbit of $\ba$ gives rise to the Fermat motive $\cM_A$,
and on the other hand, multiplications by $\hat G$-invariant monomials yields
the corresponding monomial class. These two actions are compatible.
Conversely, starting with a monomial $X_1^{v_1}X_2^{v_2}X_3^{v_3}X_4^{v_4}X_5^{v_5}$, let
$([v_1],[v_2],[v_3],[v_4],[v_5])$ be the equivalence class of monomials under multiplication 
of the $\hat G$-invariant monomials. Then there corresponds the motive
$$[[v_1]q_1,[v_2]q_2,[v_3]q_3,[v_4]q_4,[v_5]q_5]\in \frak{A}(Q).$$
Note that the constant monomial class corresponds to the weight motive $\cM_Q=[q_1,q_2,q_3,q_4,q_5]$.

This establishes the motive-monomial correspondence at the Fermat point.
\medskip

\begin{expl}{\rm
For instance for the quintic threefolds the monomial classes $(4,1,0,0,0)$ and $(3,2,0,0,0)$ combine to give only one motive $[1,1,1,3,4]$:
\medskip

\begin{center}
\begin{tabular}{|c|c|}\hline
Monomials&Motives\\\hline
$(4,1,0,0,0)$ &\\
$(0,2,1,1,1)$&\\
$(1,3,2,2,2)$ &$\in[1,1,1,3,4]$\\
$(2,4,3,3,3)$ &$\in[1,1,1,3,4]$\\
$(3,0,4,4,4)$ &\\\hline
\end{tabular}
\begin{tabular}{|c|c|}\hline
Monomials&Motives\\\hline
$(2,3,0,0,0)$& \\
$(3,4,1,1,1)$& $\in[1,1,1,3,4]$\\
$(4,0,2,2,2)$& \\
$(0,1,3,3,3)$& \\
$(1,2,4,4,4)$& $\in[1,1,1,3,4]$\\\hline
\end{tabular}
\end{center}
}
\end{expl}
\medskip

\begin{expl}\label{2030060100}{\rm
For the octic threefolds the monomial classes $(2,0,3,0,0)$ and $(6,0,1,0,0)$ together correspond to two different motives, $[4,2,2,4,4]$ and $[7,1,2,2,4]$:
\medskip
\begin{center}
\begin{tabular}{|c|c|}\hline
Monomials&Motives\\\hline
$(2,0,3,0,0)$ &\\
$(3,1,0,1,1)$&\\
$(4,2,1,2,2)$ &$\in[4,2,2,2,4]$\\
$(5,3,2,3,3)$ &$\in[7,1,2,2,4]$\\
$(6,4,3,0,0)$ &\\
$(7,5,0,1,1)$& \\
$(0,6,1,2,2)$& \\
$(1,7,2,3,3)$&$\in[7,1,4,2,2]$\\\hline
\end{tabular}
\begin{tabular}{|c|c|}\hline
Monomials&Motives\\\hline
$(6,0,1,0,0)$& \\
$(7,1,2,1,1)$& $\in[7,1,4,2,2]$\\
$(0,2,3,2,2)$&\\
$(1,3,0,3,3)$&\\
$(2,4,1,0,0)$&\\
$(3,5,2,1,1)$&$\in[7,1,4,2,2]$\\\
$(4,6,3,2,2)$&$\in[4,2,2,2,4]$\\
$(5,7,0,3,3)$&\\\hline
\end{tabular}
\end{center}
}
\end{expl}
\medskip

\begin{expl}{\rm
It should be noted that whenever a monomial class only makes a contributions at special points in 
the moduli space, such as conifold points (i.e. not at the Fermat point), none of the monomials 
in the class have no non-zero entries, and hence, as expected, no correspondence with Fermat 
motives can be made.
\medskip

\begin{center}
\begin{tabular}{|c|}\hline
$Q=(1,1,2,2,2)$\\\hline
$(4,0,3,2,1)$ \\
$(5,1,0,3,2)$\\
$(6,2,1,0,3)$ \\
$(7,3,2,1,0)$ \\
$(0,4,3,2,1)$ \\
$(1,5,0,3,2)$ \\
$(2,6,1,0,3)$ \\
$(3,7,2,1,0)$\\\hline
\end{tabular}
\begin{tabular}{|c|}\hline
$Q=(1,1,1,1,1)$\\\hline
$(4,0,3,2,1)$\\
$(0,1,4,3,2)$\\
$(1,2,0,4,3)$ \\
$(2,3,1,0,4)$ \\
$(3,4,2,1,0)$ \\\hline
\end{tabular}
\end{center}
}\end{expl}
\medskip

The above correspondence can be easily seen in the existing expressions derived for the number of points for the following cases: the one-parameter family of
the quintic Calabi--Yau threefolds by Candelas et al. [CORV00, CORV03] and for the
two-parameter family of the octic Calabi--Yau threefolds by Kadir [Ka04, Ka05]. 

\begin{expl} {\rm  
Let $m=5$ and $Q=(1,1,1,1,1)$. Then the $N_{mon}$  computed in terms of
Gauss sums formed from the Dwork character as follows.

(a) When $5|q-1$
$$N_{mon}=\frac{1}{q}\sum_{a=0}^{3}\sum_{\textbf{v}\atop v_i+a\neq 0}\lambda_{\textbf{v}}\prod_{i=1}^{5}G_{-(v_i+a)k }$$
with $k=\frac{q-1}{5}$ 

(b) When $5\nmid q-1$ 
$$N_{mon}=\frac{1}{q}G_{0}^5$$
 
The corresponding expression in terms of Jacobi sums from Weil's method is as follows:
$$\sum_{\ba\in\frakA(Q)}J(\ba)=\frac{1}{q}\sum_{\ba\in\frakA(Q)}\prod_{i=1}^5 G(\chi^{a_i}).$$
Therefore,
$$N_{motive}=N_{mon}.$$}
\end{expl}

\begin{expl}{\rm  For the two-parameter model in with $Q=(1,1,2,2,2)$ the expression for
$N_{mon}$ is more complicated. At the Fermat point, they are given as follows.

(a) When $8|q-1$
$$N_{mon}=
\frac{1}{q}\sum_{b=0}^{1}\sum_{a=0}^{3}\sum_{\textbf{v} }\lambda_{\textbf{v}}\prod_{i=1,2}G_{-(v_i+a+4b)k }
\prod_{i=3,4,5}G_{-2(v_i+a)k}$$
There is a condition on the sum over the monomial that $v_i+a+4b\neq 0$ for $i=1,2$, and $v_i+a\neq 0$ for $i=3,4,5$.

(b) When $4|q-1,\, 8 \nmid q-1$
$$N_{mon}=
\frac{1}{q}\sum_{a=0}^{3}\sum_{\textbf{v}}\lambda_{\textbf{v}}\prod_{i=1,2}G_{-\frac{(v_i+2a)}{2}l }
\prod_{i=3,4,5}G_{-(v_i+2a)l}.$$
There is a condition on the sum over the monomials that $v_i+2a\neq 0$ for $i=1,\ldots 5$.

(c) When $4\nmid q-1$
$$N_{mon}=
\frac{1}{q}\sum_{b=0}^{1}\sum_{\textbf{v}\atop v_i+4b\neq 0 }\lambda_{\textbf{v}}\prod_{i=1,2}
G_{-\frac{(v_i+4b)}{4}m }\prod_{i=3,4,5}G_{-\frac{v_i}{2}m}.$$
There is a condition on the sum over the monomial that $v_i+4b\neq 0$ for $i=1,2$, and $v_i\neq 0$ for $i=3,4,5$.
}
\end{expl}

\section{Conclusions and further problems}\label{conclusions}

In the above we had four objects to describe the topological mirror symmetry for 
Calabi--Yau hypersurfaces of dimension $3$. They are:

(a) Reflexive polytopes and their integral vertices (Toric geometry);

(b) Monomials in the graded polynomial ring, and the action of
the group of automorphisms of the coordinate ring, which gives rise to
monomial classes (Toric geometry via coordinate rings);    

(c) Picard--Fuchs differential equations (and hence periods), 
and the action of the differential operators ${\mathcal D}_i$, which gives
rise to Picard--Fuchs differential systems (Differential equations and Periods), and 

(d) Fermat motives (i.e., characters in $\frakA(Q)$ and the action of the Galois 
group $(\Z/m\Z)^{\times}$ of the $m$-th cyclotomic field $L$ over $\Q$) (Algebraic
number theory).
\smallskip

The bijection between (a) and (b) was established by Batyrev [Ba94], and 
Aspinwall--Greene and Morrison [AGM93],
the bijection between (b) and (c) follows from the work of Dwork, Katz, Griffiths (see
Cox and Katz [CK99]); the explicit bijections were constructed for the quintic Calabi--Yau
threefolds by Candelas et al. [CORV00, CORV03], and for the octic Calabi--Yau threefolds
with two deformation parameters by Kadir [Ka04, Ka05].

The point of this article is to establish a one-to-one correspondence between Fermat
motives and any of the above three objects in (a), (b) and (c) in the equivariant fashion
compatible with the various actions. 
This correspondence was established at the Fermat (the Landau--ginzburg) point in the moduli space of the 
family of Calabi--Yau threefolds.  

Though Fermat motives are defined only at the Fermat point (and also motives are algebraic
in nature), through their correspondences to monomials, they appear to contain some information 
about Calabi--Yau orbifolds with deformation parameters.  

This raises the following questions.
\smallskip

\begin{prob}
{\rm We now introduce parameters into the defining hypersurface equations of
our Calabi--Yau orbifolds. Batyrev's method produces monomials for these
Calabi--Yau hypersurfaces.  Our motives are defined at the Fermat point putting
deformation parameters equal to zero.  When deformation parameters are non-zero,
we still have monomials which correspond to Fermat motives, but acquire more monomials.  
Some monomials do arise from conifold singularities. How can one interpret
these extra monomials from the motivic point of view? 

Na\"ively, one way to start would be to group the monomials multiplied by weight 
under their transformation properties under $(\Z/m\Z)^{\times}$. For instance considering 
Example \ref{2030060100} once more and multiplying the monomials by the weight 
$(1,1,2,2,2)$ we get:
\smallskip

\begin{center}
\begin{tabular}{|c|c|}\hline
Weighted Monomials&Motives\\\hline
$(2,0,6,0,0)$ &$[2,0,6,0,0]$\\
$(6,0,2,0,0)$& \\\hline
$(3,1,0,2,2)$&$[3,1,0,2,2]$\\
$(1,3,0,6,6)$&\\
$(7,5,0,2,2)$& \\
$(5,7,0,6,6)$&\\\hline
$(4,2,2,4,4)$ &$[4,2,2,2,4]$\\
$(4,6,6,4,4)$&\\\hline
\end{tabular}
\begin{tabular}{|c|c|}\hline
Weighted Monomials&Motives\\\hline
$(5,3,4,6,6)$ &$[7,1,2,2,4]$\\
$(3,5,4,2,2)$&\\
$(7,1,4,2,2)$& \\
$(1,7,4,6,6)$&\\\hline
$(6,4,6,0,0)$ &$[6,4,6,0,0]$\\
$(2,4,2,0,0)$&\\\hline
$(0,6,2,4,4)$&$[0,6,2,4,4]$\ \\
$(0,2,6,4,4)$&\\\hline
\end{tabular}
\end{center}
\smallskip

\noindent 
Hence we obtain $4$ new ``fictitious motives'' containing zero entries, $[2,0,6,0,0]$, $[3,1,0,2,2]$, 
$[6,4,6,0,0]$ and $[0,6,2,4,4]$, in addition to the two found at the Fermat point, 
$[4,2,2,2,4]$ and $[7,1,2,2,4]$.

The  ``fictious motives'' containing zeroes may be related to lower dimensional genuine motives of Fermat 
varieties of the same degree (there may be a twisting by a character). 
For instance, we drop the component $0$ in both $[3,1,0,2,2]$ and $[0,6,2,4,4]$ 
we obtain $[3,1,2,2]$ and $[6,2,4,4]$ which may be considered as motives arising from
the Fermat surface $Z_1^8+Z_2^8+Z_3^8+Z_4^8=0\in\BP^4$ of geometric genus
$p_g=35$.  The motive $[3,1,2,2]$
has dimension $4$ and multiplicity $6$. The Hodge numbers are $h^{2,0}=2, h^{1,1}=4, h^{2,0}=2$
so that the 2nd Betti number is $B_2=8$. The motive $[6,2,4,4]$ has dimension $2$ and
multiplicity $1$. The Hodge numbers are $h^{2,0}=0, h^{1,1}=2, h^{0,2}=0$.  
However, it is not clear how these Fermat motives of lower dimensions come into the
picture of the octic Calabi--Yau threefold. 

When there is a conifold singularity, it is locally isomorphic to the quadric
$X^2+Y^2+Z^2+T^2=0$, and there should correspond a twisted Tate motive (an
extension of a Tate motive), which in turn should come from some monimial class. 

Exploring the role played by the deformation parameters will be a project in the
future. It is well known that Fermat hypersurfaces are dominated
by the product of Fermat hypersurfaces of lower dimensions of the same degree (see,
for instance Hunt and Schimmrigk [HS99]).
However, the situation is totally different if one passes onto deformations of Fermat 
hypersurfaces. Chad Schoen [Sch96] showed that the one-parameter deformation of
the quintic Calabi--Yau threefold 
$$Z_1^5+Z_2^5+Z_3^5+Z_4^5+Z_5^5-5\psi Z_1Z_2Z_3Z_4Z_5=0\in\BP^4\times\BP^1$$
is not dominated by product varieties. 
This would imply that extra monomials $(4,0,3,2,1)$ of degree $10$ and multiplicity $24$,
which corresponds to the conifold singular point $\psi=1$, are not arising 
from lower dimensional Fermat surfaces of degree $5$. At $\psi=1$, they
acquire $125$ nodes, and resolving them with $\BP^1$ yields the rigid
Calabi--Yau threefold considered by Chad Schoen [Sch86]. 
Candelas et al. [CORV03] explained the role of these extra monomials in the mirror
construction. From their calculations of the local zeta-functions, clearly we see
the appearance of the Tate motives corresponding to these extra monomials.
}
\end{prob}

\begin{prob} {\rm The study of the zeta-functions and $L$-series of Calabi--Yau orbifolds
with deformation parameters is an ongoing projct with Y. Goto and R. Kloosterman ([GKY]),
where we use a rigid cohomology theory, e.g., Monsky--Washnitzer $p$-adic
cohomology theory.  It is our hope to describe the correspondence between
monomials and motives equivariantly in terms of this $p$-adic cohomology theory 
and the other cohomology theories (e.g., \'etale, Betti). }
\end{prob}

\begin{prob} {\rm
The duality described between the two finite abelian groups
$G=G_Q$ and $\hat G$ in Section 4 may be extended further to pairs of quotients
$(X/H, \hat X/\hat H)$ where $H\subset \hat G$ and $\hat H$ is the complement
of $H$ in $\hat G$ such that $\# H\times \# \hat H=m^3$. $H$
and $\hat H$ act on $X$ and $\hat X$ respectively. 
Indeed, mirror symmetry can be extended to the mirror pairs of Calabi--Yau
orbifolds corresponding to $(H, \hat H)$ (Klemm and Theisen [KT93]), where
$H$ is normalized in some cases.

Establish motive-monomial correspondence for these Calabi--Yau threefolds. For this,
one should calculate which monomials (or equivalently, motives) are
preserved under the orbifolding operation by various different groups.  }
\end{prob}

{\bf Example 10.3.1}:
{\rm Consider the case 
$$m=5,\,Q=(1,1,1,1,1)\,:\, Y_1^5+Y_2^5+Y_3^5+Y_4^5+Y_5^5=0$$
with the groups $(G_Q,\hat G)=(\{1\},(\Z/5\Z)^3)$, and the relevant dual
pairs of subgroups. They are tabulated below where 
the number $\hat \#$ indicates a mirror partner
of $\#$.}

\begin{table}
\begin{center}
\begin{tabular}{|c||c|c|c|c|c|} \hline
 & H & generators & $h^{1,2}$ & $h^{1,1}$ & $\chi$ \\ \hline
$1$ & $\{1\}$ &    & $101$ & $1$ & $-200$ \\ \hline
$2$ & $(\Z/5\Z)$ & $(1,0,0,4,0)$ & $49$ & $5$ & $-88$ \\ \hline
$3$ & $(\Z/5\Z)$ & $(1,2,3,4,0)$ & $21$ & $1$ & $-40$ \\ \hline
$4$ & $(\Z/5\Z)^2$ & $(1,0,0,4,0),(1,2,3,4,0)$ & $21$ & $17$ & $-8$ \\ \hline
$\hat 4$ & $(\Z/5\Z)$ & $(1,2,2,0,0)$ & $17$ & $21$ & $8$ \\ \hline
$\hat 3$ & $(\Z/5\Z)^2$& $(1,2,3,4,0),(1,0,2,2,0)$ & $1$ & $21$ & $40$ \\ \hline
$\hat 2$ &$(\Z/5\Z)^2$ & $(1,0,0,4,0), (1,0,4,0,0)$ & $5$ & $49$ & $88$ \\ \hline
$\hat 1$ & $(\Z/5\Z)^3$ &$(1,0,0,4,0),(1,0,4,0,0),(1,4,0,0,0)$ &$1$& $101$ & $200$ \\ \hline
\end{tabular}
\caption{TABLE }
\end{center}
\end{table}
\smallskip

{\bf Example 10.3.2}:
{\rm Next, we consider Calabi--Yau orbifolds corresponding to
$$m=10,\,Q=(1,1,1,2,5)\,:\, Y_1^{10}+Y_2^{10}+Y_3^{10}+Y_4^5+Y_5^2=0$$
with the groups $(G_Q,\hat G)=((\Z/10\Z), (\Z/10\Z)^2)$. 
The relevant subgroups $H$ are normalized by factoring out the
group $G_Q=(\Z/10\Z)$, that is, $(\Z/2\Z)$ actually means
$(\Z/2\Z)\times(\Z/10\Z)$. }

\begin{table}
\begin{center}
\begin{tabular}{|c||c|c|c|c|c|}  \hline
 & $H$ & $\mbox{generators}$ & $h^{2,1}$ & $h^{1,1}$ & $\chi$ \\ \hline
$1$ & $\{1\}$ &                     & $145$     & $1$       & $-288$ \\ \hline
$2$ & $(\Z/2\Z)$ & $(0,5,5,0,0)$ & $99$ & $3$           & $-192$ \\ \hline
$3$ & $(\Z/2\Z)^2$ & $(0,5,5,0,0),(5,5,0,0,0)$ & $67$ & $7$ & $-120$ \\ \hline
$4$ & $(\Z/5\Z)$ & $(0,4,4,1,0)$ & $47$ & $ 11$ & $-72$ \\ \hline
$5$ & $(\Z/5\Z)$ & $(0,8,2,0,0)$ & $37$ & $13$ & $-48$ \\ \hline
$6$ & $(\Z/10\Z)$ & $(9,0,1,0,0)$ & $39$ & $15$ & $-48$ \\ \hline
$7$ & $(\Z/10\Z)$ & $(0,5,3,1,0)$ & $29$ & $17$ & $-24$ \\ \hline
$\hat 7$ & $(\Z/10\Z)$ & $(0,7,1,1,0)$ & $17$ & $29$ & $24$ \\ \hline
$\hat 6$ & $(\Z/10\Z)$ & $(0,3,3,2,0)$ & $15$ & $39$ & $48$ \\ \hline 
$\hat 5$ & $(\Z/2\Z)\times(\Z/10\Z)$ & $(5,5,0,0,0),(0,5,3,1,0)$ & $13$ & 
$37$ & $48$ \\ \hline
$\hat 4$ & $(\Z/10\Z)\times(\Z/2\Z)$ &$(9,0,1,0,0),(0,5,5,0,0)$ &$11$ & $47$ & $72$ \\ \hline
$\hat 3$ & $(\Z/\Z)^2$& $(0,8,2,0,0),(8,0,2,0,0)$ & $7$ & $67$ & $120$ 
\\ \hline
$\hat 2$ & $(\Z/10\Z)\times(\Z/5\Z)$ & $(9,0,1,0,0),(0,8,2,0,0)$ & $3$ 
& $99$ & $192$ \\ \hline
$\hat 1$& $(\Z/10\Z)^2$ & $(9,0,1,0,0),(0,9,1,0,0)$ & $1$ & $145$ & $288$ 
\\ \hline
\end{tabular}
\caption{TABLE}
\end{center}
\end{table}
\vfill
\pagebreak

\hfil{\bf Bibliography}\hfil

\hangafter1
[AGM93]
Aspinwall, P., Greene, B., and Morrison, D., {\it The monomial-divisor mirror
map}, Internat. Math. Res. Notices (1993), pp. 319--337.

\hangafter1
[Ba94] Batyrev, V., {\it Dual polyhedra and mirror symmetry for Calabi--Yau
hypersurfaces in toric varieties}, J. Alg. Geom. {\bf 3} (1994), pp. 493--535.

\hangafter1
[CLS90] Candelas, P., Lynker, M., and Schimmrigk, R., {\it Calabi--Yau
manifolds in weighted ${\BP^4}$}, Nucl. Phys. {\bf B341} (1990), pp. 383--402.

\hangafter1
[COFKM93] Candelas, P., de la Ossa, X., Font, A., Katz, S., 
and Morrison, D., {\it Mirror symmetry for two parameter models - I},
Nuclear Phys. B {\bf 416} (1994), pp. 481--538. (Also in 
[GYau94]  pp. 483--543.)

\hangafter1
[CORV00]Candelas P., de la Ossa X., and Rodriguez Villegas F, {\it Calabi--Yau manifolds over finite fields I}, arXiv:hep-th/0012233.

\hangafter1
[CORV03] Candelas, P., de la Ossa, X., and Rodriguez--Villegas, F.,
{\it Calabi--Yau manifolds over finite fields, II}, Fields Institute
Communications {\bf 38} (2003), pp. 121--157.


\hangafter1
[CK99] David A. Cox and Sheldon Katz, {\it Mirror Symmetry and Algebraic Geometry} Mathematical Surveys \& Monographs
\textbf{68}, American Mathematical Society, 1999.

\hangafter1
[D74] Deligne, P., {\it La Conjecture de Weil I}, Pub. Math. IHES {\bf 43} (1974),
pp. 273--307.

\hangafter1
[DM82] Deligne, P., and Milne, J.S., {\it Tannakian categories, Hodge cycles, motives
and Shimura varieties}, LNM {\bf 900}, Springer--Verlag, 1982. 

\hangafter1
[Dol82] Dolgachev, I., {\it Weighted projective spaces} in ``Group Actions and Vector Fields''
LNM {\bf 956} (1982), pp. 24--71. Springer--Verlag.

\hangafter1
[Dw60] Dwork, B., {\it On the rationality of the zeta-function of an algebraic
variety}, Amer. J. Math. {\bf 82} (1960), pp. 631--648.

\hangafter1
[G06] Goto, Y., {\it  On the height of Calabi--Yau threefolds}, in preparation.

\hangafter1
[GKY] Goto, Y., Kloosterman, R., and Yui, N., {\it Zeta-functions and $L$-series of certain
$K3$-fibered Calabi--Yau threefolds}, in preparation.

\hangafter1
[GY95] Gouv\^ea, F., and Yui, N., {\it Arithmetic of Diagonal Hypersurfaces
over Finite Fields}, London Math. Soc. Lecture Notes Series {\bf 209},
Cambridge University Press 1995.

\hangafter1
[GP90] Greene, B., and Plesser, M.R., {\it Duality in Calabi--Yau
moduli spaces}, Nuclear Physics B {\bf 338} (1990), no. 1, pp. 15--37.
 
\hangafter1
[GRY91] Greene, B. R., Roan, S.-S. and Yau, S.-T., {\it Geometric singularities
and spectra of Landau--Ginzburg models}, Commun. Math. Phys. {\bf 142} 
(1991), pp.  245--259.

\hangafter1
[GYau96] Greene, R., and Yau, S.-T., {\it Mirror Symmetry II},
AMS/IP Studies in Advanced Mathematics, {\bf 1},
American Mathematical Society and International Press 1996.

\hangafter1
[HLY95] Hosono, S., Lian, B.-H., and Yau, S.-T., {\it Picard--Fuchs
differential operators for Calabi--Yau hypersurfaces with
$h^{1,1}\leq 3$}--Optional appendix to {\it GKZ--generalized hypergeometic
system in mirror symmetry for Calabi--Yau hypersurfaces}, Commun. Math. Phys.
{\bf 182} (1996), pp. 525--578.

\hangafter1
[HS99] Hunt, B., and Schimmrigk, R., {\it $K3$-fibered Calabi--Yau threefolds. I. The
twist map}, Internat. J. Math. {\bf 10} (1999), no.7, pp. 833--869.

\hangafter1
[Ka04] Kadir, S., {\it The Arithmetic of Calabi--Yau Manifolds and Mirror Symmetry}, arxiv: hep-th/0409202.

\hangafter1
[Ka05] Kadir, S., {\it Arithmetic mirror symmetry for a two-parameter
family of Calabi--Yau manifolds}, in ``Mirror Symmetry V'', AMS/IP Studies
in Advanced Mathematics, to appear.

\hangafter1
[KS] Kreuzer, M., and Skarke, H., {\it Calabi--Yau data}, 

http://hep.itp.tuwien.ac.at/kreuzer/CY/CYcy.html.

\hangafter1
[KT93] Klemm, A., and Theisen, S., {\it Considerations of one-modulus
Calabi--Yau compactifications: Picard-Fuchs equations, K\"ahler potentials
and mirror maps}, Nucl. Physics {\bf B 389} (1993), pp. 153--180.

\hangafter1
[Ma70] Manin, Yu. I., {\it Correspondences, motives and monoidal transformations},
Math. USSR-Sb. {\bf 77} (1970), pp. 475--507.

\hangafter1
[Mor97] Morrison, D. M., {\it Mathematical aspects of mirror symmetry}, in
``Complex Algebraic Geometry (Park City, Utah, 1993), AMS/Park City Math. Ser. {\bf 3}
(1997), pp. 265--327. 

\hangafter1
[Mo98] Moore, G., {\it Arithmetic and Attractors}, hep-th/9807087.

\hangafter1
[Na75] Nakamura, I., {\it Complex parallelisable manifolds and their small
deformations}, J. Diff. Geom. {\bf 10} (1975), pp. 85--112.

\hangafter1
[Ro90] Roan, S.-S., {\it On Calabi--Yau orbifolds in weighted projective
spaces}, International J. Math., {\bf 1}, No.2 (1990), pp. 211--232.

\hangafter1
[Ro91] Roan, S.-S., {\it The mirror of Calabi--Yau orbifold},
Internat. J. Math., {\bf 2}, No.4 (1991), 439--455.

\hangafter1
[Ro94] Roan, S.-S., {\it On $c_1=0$ resolution of quotient singularities},
Internat. J. Math., {\bf 4} (1994), pp. 523--536.

\hangafter1
[Sch86] Schoen, C., {\it On the geometry of a special determinantal hypersurface
associated to the Mumford--Horrocks vector bundle}, J. Reine Angew. Math. {\bf 364}
(1986), pp. 85--111.
 
\hangafter1
[Sch96] Schoen, C., {\it Varieties dominated by product varieties},
International J. Math. {\bf 7}, No. 4 (1996), pp. 541--571.

\hangafter1
[Sh87]
Shioda, T., {\it Some observations on Jacobi sums}, Adv. Stud. Pure Math. {\bf 12}
(1987), pp. 119--135.

\hangafter1
[So84] Soul\'e, C., {\it Groupes de Chow et $K$-th\'eorie de variet\'es sur un corps fini},
Math. Ann. {\bf 268} (1984), pp. 317--345.

\hangafter1
[We49] Weil, A., {\it Number of solutions of equations over finite fields},
Bull. Amer. Math. Soc. {\bf 55} (1949), pp. 497--508.

\hangafter1
[Y05] Yui, N., {\it The $L$-series of Calabi--Yau orbifolds of CM type} with
Appendix by Y. Goto, in {\it Mirror Symmetry V}, AMS/IP Studies in Advanced Mathematics, to
appear.
\vfill
\end{document}